\documentclass[leqno,14pts]{amsart}
\numberwithin{equation}{section}
\usepackage{pdfsync}
\usepackage{amsmath,amsthm,amsfonts,amssymb,amsxtra,graphicx}

\usepackage{hyperref}
\usepackage[alphabetic,initials,nobysame]{amsrefs}

\hypersetup{colorlinks=true,citecolor=blue,linkcolor=blue,urlcolor=blue,pdfstartview=FitH,
pdfauthor=Valentin Blomer and Philippe Michel,pdftitle=Sup-norms of eigenfunctions on arithmetic ellipsoids}

\theoremstyle{plain}

\newtheorem{thm}{Theorem}  
\newtheorem{lemma}{Lemma}

\newtheorem*{thm*}{Theorem}

\newtheorem*{prop*}{Proposition}
\newtheorem*{lem*}{Lemma}

\theoremstyle{definition}

\newtheorem*{example*}{Example}

\theoremstyle{definition}

\newtheorem{rem}{Remark}[section]

\newtheorem*{rem*}{Remark}

\newcommand{\x}{{\mathbf{x}}}
\newcommand{\nr}{{\mathrm {nr}}}

\newcommand{\Cl}{{\mathrm{Cl}}}

\newcommand{\mfI}{{\mathfrak{I}}}
\newcommand{\mfn}{{\mathfrak{n}}}
\newcommand{\mfp}{{\mathfrak{p}}}
\newcommand{\GL}{{\mathrm{GL}}}

\newcommand{\Nr}{\nr}
\newcommand{\ov}{{\rm nr}}
\newcommand{\vphi}{\varphi}
\newcommand{\tphi}{{\tilde{\varphi}}}
\newcommand{\tpi}{{\tilde{\pi}}}
\newcommand{\tr}{{\rm tr}}

\newcommand{\mcL}{{\mathcal L}}

\newcommand{\B}{\rmB}

\newcommand{\rmG}{\mathrm{G}}

\newcommand{\GA}{\rmG(\Aa)}
\newcommand{\GAf}{\rmG(\Af)}

\newcommand{\tGA}{\tG(\Aa)}
\newcommand{\tG}{\widetilde{\rmG}}

\newcommand{\PB}{{\mathrm{PB}^\times}}
\newcommand{\Bu}{\mathrm{B}^{1}}
\newcommand{\PBA}{{\PB(\Aa)}}

\newcommand{\PBAf}{{\PB(\Af)}}

\newcommand{\BA}{{\rmB(\Aa)}}

\newcommand{\Bt}{\rmB^{\times}}

\newcommand{\BtA}{{\rmB^{\times}(\Aa)}}
\newcommand{\BtAf}{\Bt(\Af)}

\newcommand{\rmK}{{\mathrm{K}}}

\newcommand{\eps}{\varepsilon}

 \newcommand{\ra}{\rightarrow}
 
\def\peter#1{\langle #1\rangle}
\def\ov#1{\overline{#1}}

 \DeclareFontFamily{OT1}{rsfs}{}
\DeclareFontShape{OT1}{rsfs}{n}{it}{<-> rsfs10}{}
\DeclareMathAlphabet{\mathscr}{OT1}{rsfs}{n}{it}

 \newcommand{\beq}{\begin{displaymath}}
\newcommand{\eeq}{\end{displaymath}}
\newcommand{\be}{\begin{equation}}
\newcommand{\ee}{\end{equation}}

\newcommand{\refs}{\eqref}

\newcommand{\order}{\mathscr{O}}

\newcommand{\mcH}{\mathcal{H}}
\newcommand{\vol}{{\mathrm{vol}}}
\newcommand{\rmB}{{\mathrm{B}}}
\newcommand{\rmZ}{{\mathrm{Z}}}
\newcommand{\Oo}{\order}

\newcommand{\whOt}{\widehat{\order}^{\times}}
\newcommand{\whO}{\widehat{\order}}
\newcommand{\Cc}{{\mathbb{C}}}
\newcommand{\Hh}{{\mathbb{H}}}

\newcommand{\bfm}{{\mathbf{m}}}
\newcommand{\bfl}{{\mathbf{l}}}

\newcommand{\Zz}{{\mathbb{Z}}}

\newcommand{\Rr}{{\mathbb{R}}}
\newcommand{\Qq}{{\mathbb{Q}}}

\newcommand{\Aa}{{\mathbb{A}}}
\newcommand{\At}{{\mathbb{A}}^{\times}}

\newcommand{\Af}{{\mathbb{A}}_{f}}

\DeclareMathOperator{\SL}{SL}
\newcommand{\SO}{\mathrm{SO}}
\newcommand{\TR}{\mathrm{T}(\Rr)}

\newcommand{\SU}{\mathrm{SU}_2(\Cc)}
\newcommand{\disc}{\mathrm{disc}}
\newcommand{\rdisc}{\mathrm{disc}^{\ast}}
\newcommand{\bash}{\backslash}
\newcommand{\N}{{\nr_F}}
\newcommand{\NE}{{\nr_E}}
\newcommand{\resprod}{\mathop{{\prod}^{\mathbf{'}}}\limits}

\begin{document}

\title{Hybrid bounds for automorphic forms on ellipsoids over number fields}

\begin{abstract} We prove upper bounds for Hecke-Laplace eigenfunctions on certain Riemannian mani\-folds $X$ of arithmetic type, uniformly in the eigenvalue and the volume of the manifold. The manifolds under consideration are $d$-fold products of $2$-spheres or $3$-spheres, realized as adelic quotients of quaternion algebras over totally real number fields. In the volume aspect we prove a (``Weyl-type") saving of $\text{vol}(X)^{-1/6+\varepsilon}$. 
\end{abstract}

\author{Valentin Blomer}
\address{Mathematisches Institut, Universit\"at G\"ottingen, Bunsenstr. 3-5, 37073 G\"ottingen} \email{blomer@uni-math.gwdg.de}
\author{Philippe Michel\\\today}
\address{EPFL/SB/IMB/TAN, Station 8, CH-1015 Lausanne, Switzerland}
\email{philippe.michel@epfl.ch}
\thanks{V.\ B.\ is supported by     a Volkswagen Lichtenberg fellowship and an ERC starting grant. Ph.\ M.\ is partially supported by the ERC advanced research grant n. 228304 and the  SNF grant 200021-12529.}


\maketitle
\tableofcontents
\section{Introduction}

Given a Riemannian manifold $X$, it is a classical problem to give a pointwise upper bound  for a $L^2$-normalized  Laplace eigenfunction $\phi$ in terms of the Laplace eigenvalue $\lambda$ and/or properties of $X$. It can be either seen as a rough measure of the non-concentration of the mass of $\phi$ or as a degenerate restriction problem (where the cycle is reduced to a single point). If $X$ is compact, generic methods give the bound \cite{Sa}
\begin{equation}
\label{genericbound}
\| \phi \|_{\infty}\ll_X (1+|\lambda|)^{(\dim X-1)/4},
\end{equation}
 and one seeks improvements over this bound. As pointed out in \cite{Sa}, the sup-norm problem is also closely tied to the multiplicity problem: If $V_{\lambda}$ denotes the eigenspace of the eigenvalue $\lambda$, generic methods show  that 
\begin{equation}\label{trivial}
 \dim V_{\lambda} \leq \text{vol}(X) \, \sup_{\substack{\phi \in V_{\lambda}\\ \| \phi \|_2 = 1}} \| \phi \|_{\infty}^2.
\end{equation} 
In other words, high multiplicity of eigenvalues may be an obstruction to   sup-norm bounds better than \eqref{genericbound}, at least for general eigenfunctions $\phi$. For instance, in the case of the sphere $X = S^2$, the dimension of the eigenvalue $\lambda = k(k+1)$ ($k \in \Bbb{N}_0$) is known to be $2k+1$. Hence by \eqref{trivial} the best possible sup-norm bound we can hope for is of order $(1+k)^{1/2} \asymp (1 + \lambda)^{1/4}$ which is realized by  the $L^2$-normalized zonal spherical function  $\sqrt{(2k+1)/(4\pi)} p_k(\cos \theta) $ (where $\theta$ is the polar angle and $p_k$ is the $k$-th Legendre polynomial). This situation is in sharp contrast to negatively curved Riemann surfaces where the sup-norms of eigenfunctions (and in particular multiplicities of eigenvalues) are believed to be small ($\ll_\eps (1+|\lambda|)^\eps$ for any $\eps>0$). However, even in that case, all that one knows in this generality about the multiplicity is an extremely modest (by a factor $\log(2+|\lambda|)$) improvement of Berard  (see \cite{Schur}). 

An obvious way to try to resolve the multiplicity issue is to exploit extra symmetries and to require that $\phi$ is an eigenfunction of additional operators, commuting with $\Delta$. For instance, if $X$ is locally symmetric of rank $>1$ and $\phi$ is an eigenfunction of the full algebra of invariant differential operators the bound improves \eqref{genericbound} to $\ll (1+|\lambda|)^{(\dim X-\mathrm{rank} X)/4}$ \cite{Sa}.

 This does not help much for surfaces, but if these are of "arithmetic type" and so endowed with the action of a suitable algebra of commuting Hecke operators, one may then consider instead {\em joint Hecke-Laplace eigenfunctions}, and some significant saving is possible. Indeed when $X := \Gamma \backslash \Bbb{H}$ is a modular or a Shimura curve  ($\Bbb{H} \simeq \SL_2(\Bbb{R})/\SO_2$ the hyperbolic plane and $\Gamma \subseteq SL_2(\Bbb{R})$ an arithmetic lattice),  the bound \refs{genericbound} was improved by a  power of $1+|\lambda|$ in the groundbreaking work of Iwaniec and Sarnak  \cite{IS}. Later, a similar result was obtained by Vanderkam  for $\Hh$ replaced by the $2$-sphere $S^2 \simeq \SO_3(\Rr)/\SO_2(\Rr)$ \cite{VdK}. In these cases it follows from the multiplicity one theorem and the Jacquet-Langlands correspondence that the dimension of a Hecke-Laplace eigenspace is bounded by $\ll_\eps \vol(X)^\eps$ for any $\eps>0$.

\subsection{Bounds on $2$-dimensional ellipsoids} Our previous work \cite{BM1} 
dealt with a family of varieties $X=X(L,q)$ associated to pairs $(L,q)$ for $q$ a definite quadratic form on a three dimensional $\Qq$-vector space $V$ and $L$ a suitable lattice in $V$; $X$ was then a finite union of (quotients) of $2$-spheres indexed by a set of representatives of genus classes of $L$. 

The present paper extends \cite{BM1} in two further directions: on the one hand, by a new treatment of the amplifier we improve significantly the main result in \cite{BM1}; on the other hand, we extend the argument to varieties attached to quadratic lattices $L\subset V$ for $(V,q)$ a totally definite ternary quadratic space defined over some fixed, totally real number field $F$ of degree $d$ over $\Bbb{Q}$; the corresponding variety $X(L,V)$ is then a finite union of $d$-fold products of $2$-spheres. We stress that this extension to number fields is not solely for the sake of generality: in the next subsection, we use these results to study similar problems for varieties associated to {\em quaternary quadratic spaces}.

Our main results are proven under some additional assumption on the "shape" of $L$ which is better expressed in terms of quaternions (we refer to \S \ref{quaternion} for the notations related to quaternions). Given  a totally real  field $ F/\Qq$ of degree $d$, a  totally definite ternary quadratic space $(V, q)$ over $F$ and a lattice $L\subset V$, there exists (cf.\ \cite[Chap.\ 1]{kitaoka}) $\lambda\in F^\times$ and a (unique up to isomorphism) totally definite quaternion algebra $\rmB$ over $F$ (i.e.\ at all real places of $F$, $\B$ is isomorphic to the  Hamilton quaternions) such that $(V,q)$ is isometric to the ternary quadratic space $(\B^0,\lambda \nr)$; here $\B^0$ denotes the space of trace $0$ quaternions and $\nr$ the reduced norm of $\B$. Choosing such an isometry, we therefore identify $V$ with $\B^0(F)$ and $L$ with a certain sublattice of $\B^0(F)$. We assume from now on that $$L=\Oo^0:=\B^0(F)\cap\Oo$$
for  an Eichler order $\Oo$ of $\B(F)$. 

To this situation is associated a finite disjoint union of quotients of $d$-fold products of $2$-spheres
\begin{equation}\label{defX2}
 X^{(2)}(\Oo):=\bigsqcup_{i\in I}X_{i},\quad  X_{i}:=\Gamma_{i}\bash (S^{2})^d,
\end{equation}
for $\Gamma_{i}< \SO_{3}(\Rr)^d$ some finite subgroup of order bounded in terms of $d$ only. 
This union corresponds to a double adelic quotient whose definition is given in \eqref{defX2adelic}.
 
The  quotients $X_{i}$ are called the components of $X^{(2)}(\Oo)$ and their indexing set is 
 the set of classes of left-$\Oo$ ideals in $\B$ (whose cardinality is the class number of $\Oo$). Consider the restriction to $S^2$ of the Euclidean metric on $\Rr^3$; this is an $\SO_3(\Rr)$-invariant Riemannian metric on $S^2$ which induces a volume form and a $d$-tuple of Laplace operators $\Delta=(\Delta_1,\dots,\Delta_d)$ on $(S^2)^d$ that descend to $X^{(2)}(\Oo)$. One has (see \S \ref{adelicinterp} below)
\begin{equation}\label{defV2}
V_2:= \vol(X^{(2)}(\Oo)))=|\Nr_{F/\Bbb{Q}}\disc(\Oo)|^{1/2+o_F(1)}
\end{equation}
(here the $\disc(\Oo)$ refers to the discriminant of the quaternary space $(\Oo,\Nr)$).

We are interested in obtaining non-trivial bounds for  the $L^\infty$-norm of an $L^2$-normalized $\Delta$-eigenfunction $\vphi$ on $X^{(2)}(\Oo)$ in terms of the eigenvalues $\lambda=(\lambda_1,\dots,\lambda_d)$ of $\Delta$ and of the total volume $\vol(X^{(2)}(\Oo))$.  The trivial bound in this case is (see \cite{Sa} for a general result)
\begin{equation}\label{deflambda} \|\vphi\|_\infty \ll |\lambda|^{1/4},\ \text{with\footnotemark} \,\, |\lambda|:=\prod_{j=1\dots d}(1+|\lambda_j|).\end{equation}\footnotetext{Formally speaking, the symbol $|\lambda|$  should be thought of as a single quantity, not the absolute value of some real number.}

Our objective is to improve over this bound \emph{simultaneously} in the $\lambda$ and the volume aspect; such non-trivial bounds are called ``hybrid''. In this generality this is hopeless: the previous bound is indeed sharp both in the volume and in the $\lambda$-aspect. The possibility of constructing Laplace eigenfunctions with large sup-norm comes from the fact that $\Delta$-eigenfunctions have very large multiplicities (roughly $\approx V_2|\lambda|^{1/2})$. As explained above, a way to resolve this issue is to require $\vphi$ to be also an eigenfunction of a family of ``Hecke" operators, indexed by the complement of a finite, fixed subset of the prime ideals of $F$, $\{T_\mfp\}_{\mfp\nmid \disc(\Oo)}$. 
The Hecke operators $\{T_\mfp,\ \mfp\nmid \disc(\Oo)\}$ together with $\Delta$ generate a commutative algebra of self-adjoint operators on $L^2(X^{(2)}(\Oo))$; in particular this space admits an orthonormal basis made of Laplace-Hecke eigenfunctions.
 
\begin{thm}\label{mainthm} Let $\Oo$ be an Eichler order in a totally definite quaternion algebra $\B$ over $F$, and let $\vphi$ be an $L^{2}$-normalized Hecke-Laplace  eigenfunction on $X^{(2)}(\Oo)$.  
Then one has with the notation as in \eqref{deflambda}, 
\begin{equation}\label{hybrid}
\|\vphi\|_\infty \ll |\lambda|^{\frac{1}{4}} (V_2|\lambda|^{\frac{1}{2}})^{-\frac{1}{20}}.
\end{equation}
Individually, we obtain the following bounds in the $\lambda$ and in the volume aspect
\begin{equation}\label{indiv}
  \|\vphi\|_\infty \ll_{\varepsilon}  |\lambda|^{\frac{1}{4}} V_2^{-\frac{1}{6}+\varepsilon},\ \|\vphi\|_\infty \ll  |\lambda|^{\frac{1}{4} -\frac{1}{27}}.
 \end{equation} 
\end{thm}

The bound \eqref{hybrid} is obtained by interpolation between the two bounds in \eqref{indiv}. We emphasize that as in \cite{BM1} these estimates are uniform in $\B$ and $\Oo$, but we regard the number field $F$ as fixed. For the rest of the paper all implied constants may depend on $F$ as well as on a small real number $\varepsilon$ where appropriate. All other dependencies will be indicated. \\


The first non-trivial bound of this sort was obtained  by Iwaniec and Sarnak \cite{IS} for $F=\Qq$, $\B$ indefinite,  and a fixed order $\Oo$. For $\Oo$ varying (of square-free level) a bound simultaneously non-trivial in $\vol(X^{(2)}(\Oo))$ and in  $|\lambda|$ was obtained by the first named author and R. Holowinsky  \cite{BH}. This result was extended by Templier \cite{Te} to the case of a totally real number field and for $\B$ indefinite at one archimedean place. 
 In the definite case, the first non-trivial result we are aware of is due to Vanderkam \cite{VdK}: for $F=\Qq$, $\B$ the Hamilton quaternions and $\Oo$ the maximal order, he obtained
\begin{equation}\label{VanderKam}
  \|\vphi\|_\infty \ll |\lambda|^{\frac{1}{4} - \frac{1}{24}+\varepsilon}.
\end{equation}   
Unaware of his work, we proved in \cite{BM1} a hybrid bound for general $\B$ and any Eichler order $\Oo$ of the shape
\begin{displaymath}
\|\vphi\|_\infty \ll |\lambda|^{\frac{1}{4}} (V_2|\lambda|^{\frac{1}{2}})^{-\frac{1}{60}+\varepsilon},
\end{displaymath}
which was an interpolation between the individual bounds 
$$\|\vphi\|_\infty \ll |\lambda|^{1/4} V_2^{-1/12+\varepsilon}\hbox{ and }\|\vphi\|_\infty \ll  |\lambda|^{11/48+\varepsilon} V_2^{1/12+\varepsilon}.$$Our present result  \eqref{indiv} is stronger  in both aspects.   
The improvement in the volume aspect  comes from a new way to deal with the amplifier (occurring from the amplification method) which may be of general interest. In the $\lambda$ aspect, the improvement comes from the use of Vanderkam's method. Our bound in \eqref{indiv}, however,  is marginally weaker than \eqref{VanderKam} because of some technical obstacles in the number field case. \\

We remark that the strongest conceivable result in the situation of Theorem \ref{mainthm} is $$\|\vphi\|_{\infty} \ll |\lambda|^{1/4} (V_2|\lambda|^{1/2})^{-1/2+\varepsilon}.$$ It seems reasonable to conjecture that in the presently considered case of a compact manifold this bound reflects the reality, although some care has to be taken as N.\ Templier \cite{Tem1} has recently disproved a similar conjecture in the non-compact case. 
 
 In any case, in the level aspect, we arrive at least at $33\%$ of the true bound which is similar to Weyl's bound vs.\ the Lindel\"of Hypothesis for Riemann's zeta function. Gergely Harcos and Nicolas Templier kindly informed us that 
for the (indefinite) discriminant quadratic form $b^2-4ac$ over $\Qq$, they obtained in \cite{HT} the   same exponent $-1/6$. This convergence of exponents obtained independently and in fairly different contexts makes it therefore likely that this result will be hard to improve with the present technology.

\subsection{Application to $3$-dimensional ellipsoids}
We illustrate the extension of \cite{BM1} to general totally real number fields by providing non-trivial sup-norm bounds for Hecke-Laplace eigenfunctions on manifolds $X$ that are finite unions of $d$-fold products of \emph{$3$-spheres} $S^3 = \SO_4/\SO_3$ (i.e.\ bounds for automorphic forms of orthogonal groups in $4$ variables). The main point here is that there is a close relationship between automorphic forms on orthogonal groups in $4$ variables over $F$ and automorphic forms on orthogonal groups in $3$-variables over a suitable (possibly split) quadratic extension  $E$ of $F$ (an extension of the well known fact that $\SO_4(\Rr)$ is a double cover of $\SO_3(\Rr) \times \SO_3(\Rr)$.)

As above these manifolds are better described in terms of quaternion algebras: recall (see \S \ref{sec61}) that to any non-degenerate quadratic space $(V,Q)$ over $F$ with discriminant $\Delta$, there is canonically associated a quaternion $F$-algebra $\B$, a quadratic etale $F$-algebra $E$ ($F\times F$ if $\Delta$ is a square in $F^\times$ and $F(\sqrt\Delta)$ otherwise) and a four dimensional vector space $$\B'\subset \B_E=\B\otimes_ F E$$ such that $(V,Q)$ is similar\footnote{Even isometric if $F$ is totally real and $q$ is positive at every archimedean place, by Eichler's norm theorem} to $(B', \Nr)$. 

Given an Eichler order $\Oo\subset\B$, we associate to it the integral quadratic lattice $$(\Oo', \Nr)\hbox{ for } \Oo'=\B'\cap \Oo\otimes_{{\mathcal{O}_F}}\mathcal{O}_E$$
and a finite disjoint union of (quotients of) $d = [F : \Bbb{Q}]$ products of $3$-spheres (cf.\ \eqref{defX3})
$$
 X^{(3)}(\Oo)=\bigsqcup_{i\in I}X_{i},\quad  X_{i}:=\Gamma_{i}\bash (S^{3})^d,
$$
where $\Gamma_{i}< \SO_{4}(\Rr)^d$ is finite and of order bounded in terms of $d$ only. Again the indexing set is closely related to the set of genus classes of the quaternary quadratic lattice $(\Oo', \Nr)$ and if one equips $S^3$ with the restriction of the Euclidean metric on $\Rr^4$, the total volume of $X^{(3)}(\Oo)$ satisfies
\begin{displaymath}
  V_3:= \vol(X^{(3)}(\Oo))= |\Nr_{F}(\disc(\Oo))|^{1+o_E(1)}.
  \end{displaymath}
Similarly as above, $X^{(3)}(\Oo)$ is endowed with a commutative algebra of Hecke operators commuting with the corresponding Laplace operator. For  an  $L^2$-normalized Hecke-Laplace eigenfunction $\vphi$ on $X^{(3)}(\Oo)$ the trivial bound for its sup-norm is  $$\|\vphi\|_{\infty} \ll |\lambda|^{1/2},$$
(uniformly in the volume $V_3$) and we obtain here an  improvement in the volume aspect:
\begin{thm}\label{thmso4} In the situation described above, one has
\begin{equation}\label{so4}
\| \vphi \|_{\infty} \ll_{E, \varepsilon} |\lambda|^{1/2} V_3^{-1/6+\varepsilon}
\end{equation}
 for any $\eps>0$. 
\end{thm}

\medskip

The present bound is a direct application of the arguments of the proof of Theorem \ref{mainthm}; yet it seems to be the first instance of a non-trivial arithmetic (i.e.\ in the level aspect) sup-norm bound for a manifold which does not factor into surfaces. Again we obtain the same (``Weyl-type") quality in the exponent. Several extensions are possible:

\begin{enumerate}
\item We have considered here only the volume aspect. The diophantine counting Lemma    \ref{VdKlemma} of \S \ref{sec:quadratic} of this paper would yield quite directly some non-trivial hybrid bounds for some $\SO_Q$-automorphic forms, namely those, which at each archimedean place of $F$ correspond (via the identification $\SO_Q(F_\sigma)\simeq \SO_4(\Rr)$) to pure weight vectors with respect to the action of the maximal torus $\SO_2(\Rr)\times\SO_2(\Rr)<\SO_4(\Rr)$. Laplace-Hecke eigenfunctions on $3$-dimensional ellipsoids on the other hand, correspond to $\SO_3(\Rr)$-invariant vectors; these are potentially {\em long} linear combinations of pure weight vectors, and  Lemma 
 \ref{VdKlemma} in its present form is not sufficient to obtain hybrid bounds for such functions.
 
\item  The present bound depends on the quadratic extension $E$. Making it explicit and non-trivial in this aspect   requires a more precise description of the
local structure of the  quaternary quadratic lattices considered at the places where $E$ is ramified and versions of the counting Lemmata \ref{lemma2} -- \ref{VdKlemma} taking this aspect into account. Observe that in the present case,  the amplification method  does not a priori require that $E$ splits at many small places (as is the case in \cite{DFIDuke} or \cite[\S 7]{Ven}), for the group $\SO_Q(F_v)$ has rank at least $1$ for almost all places of $F$ (the places at which $\B$ is unramified).

\end{enumerate}

\subsection{Organization of the paper and concluding remarks} In the next section we introduce general notations and describe how the problem translates in the adelic setting. Section \ref{reduction} discusses reduction theory for totally definite quadratic forms over totally real number fields, and we discuss general results about the representation of algebraic integers by such quadratic forms in Section \ref{sec:quadratic}. 
In Section \ref{pretrace} we apply the pretrace formula and the amplification method in a by now standard way and reduce the problem of bounding the sup norm of Hecke-Laplace eigenforms to the diophantine problems of the previous section, that is,   bounding representation numbers of quadratic forms of  large discriminant of $F$-integral vectors that are almost parallel or almost orthogonal to a given vector. 
The first bound in \eqref{indiv} and the bound \eqref{so4} follow only from Lemma \ref{lemma2} which is at least in principle not much more than a generalized Lipschitz principle. The second bound in \eqref{indiv} is more complicated and requires Lemmata \ref{lemma3} -- \ref{VdKlemma}. There are at least two sources of improvement compared to the analysis in \cite{BM1}:  in the present paper we use explicitly the fact that the considered quadratic forms are associated to an order in a quaternion algebra and in particular represent $1$. Moreover, we exploit the average over the amplifier and treat the quadratic part of variable $\ell$ in the amplifier essentially as a new variable of the quadratic form. \\

Finally we would like to thank the referee   for his unusually careful reading of the manuscript and a long list of constructive suggestions that greatly  improved the presentation.

\section{Preliminaries}
Let $F/\Qq$ be a totally real number field of degree $d$, $\mathcal{O}_F$  its ring of integers, $U=\mathcal{O}_F^{\times}$ its group of units and $U^+$ the subgroup of totally positive units. For  a place $v$ of $F$, we denote by $F_v$ the associated local field. 
A typical real place of $F$ will be denoted as an embedding 
\begin{displaymath}
  \sigma: F \hookrightarrow \Bbb{R}
 \end{displaymath}
 and the list of real places will be denoted by $\sigma_1, \ldots, \sigma_d$; for $x\in F$, we write $x^\sigma=\sigma(x)\in\Rr$ for the corresponding conjugate.  

We denote by $\Aa=\resprod_{v}F_v$, $\Af=\resprod_{v<\infty}F_v$, $F_\infty=\prod_{\sigma} F_\sigma$ the $F$-algebras of ad\`eles, finite ad\`eles and archimedean components of $F$. We denote the   norm on $F/\Qq$ by   $\N$   and use the same notation for the natural extension of the norm to the $F$-ideals or to various $F$-algebras related to $F$ ($F_v$, $F_\infty$, $\Aa$ etc.). To ease notations, for $\mfI\subset F$ an $F$-ideal or $x_f\in\Af$ a finite $F$-id\`ele we will freely identify $\N(\mfI)$ or $\N(x_f)$ with the positive generator of its underlying  $\Zz$-ideal in $\Qq$, so that if necessary the expression $\N(\mfI)^{\sqrt \pi}$ (say) is well-defined.\\

\subsection{Ternary quadratic spaces and  quaternion algebras}\label{quaternion}  We recall some facts about quaternion algebras, see e.g.\ \cite{Vig} or \cite{BM1} for more details.  Let $\B$ be a totally definite quaternion algebra defined over $F$. We denote its canonical involution by $z\mapsto z^*$ and its reduced trace and reduced norm by 
$$\tr: z\mapsto z+z^*,\ \nr: z\mapsto zz^*.$$ We denote the trace-$0$ quaternions and trace-$0$ quaternions of norm $1$ by $\B^0$ and $\B^{0,1}$ respectively (considered as algebraic varieties over $F$). The spaces $(\B^0,\nr)$, $(\B,\nr)$ are  quadratic $F$-spaces whose associated inner product is denoted by
$$\peter{z_1,z_2}_\B=\frac12\tr(z_1 z_2^*),\ \peter{z,z}_\B=\nr(z).$$

 We denote by $\Bt$ the group of units, $\rmZ$ its center (the subgroup of scalars), $\B^1$ the subgroup of quaternions of reduced norm $1$ and by $\PB= \rmZ\bash\Bt$ the projective quaternions. All these are considered as $F$-algebraic groups in the evident way.  We write $\BA$, $\BtA$, $\Bu(\Aa)$, $\BtAf$, $\ldots$, $\B(F_v),\ \Bt(F_v)$ etc.\ for the sets of rational points of these varieties over the corresponding $F$-algebras.
 
The conjugation action of the group of units $\Bt$ on the ternary quadratic space $(\B^0,\nr)$ is isometric (i.e.\ preserves the norm form) and the map
$$g\in\Bt \mapsto \rho_g:\begin{array}{ccc}\B^0 & \ra & \B^0 \\x & \mapsto & gxg^{-1}\end{array}$$
is an isomorphism of $F$-algebraic groups 
$$\PB=\rmZ\bash\Bt\simeq \SO(\B^0),\quad  \Bu\simeq \mathrm{Spin}(\B^0)$$
where $\mathrm{Spin}(\B^0)$ denotes the spin group (the simply connected covering group of $\SO(\B^0)$).

A place $v$ is called ramified if $\rmB_v := \rmB \otimes_{F} F_v$ is a division algebra, and non-ramified otherwise; in the former case, $\rmB_{v}$ is the unique (up to isomorphism) quaternion division algebra over $F_v$; in the latter $\rmB_v \cong \text{Mat}(2, F_v)$ in which case the reduced norm and reduced trace are given by the usual determinant and trace for matrices. 
Since $\B$ is totally positive, all archimedean places are ramified.   The {\em reduced discriminant} $D_{\rmB}$ of $\rmB$ is the product of the finite ramified prime ideals. 

A lattice or ideal $\mfI\subset \rmB$ is an $\mathcal{O}_F$-module of maximal rank $4$. An order $\Oo$ is a subring of $\B$ which is also lattice. The left   order $\Oo_{l}(\mfI)$  of a lattice $\mfI$ is the set
$\Oo_{l}(\mfI)=\{\gamma\in \B,\ \gamma\mfI\subset\mfI\}$.  Given an order $\Oo$, a left $\Oo$-ideal is defined as a lattice $\mfI$ such that $\Oo_{l}(\mfI)=\Oo$. Two left $\Oo$-ideals $\mfI,\mfI'$ are called (right-)equivalent if there exists $\gamma\in\Bt(F)$ such that $\mfI'=\mfI \gamma$. The set of such equivalence classes is denoted $\Cl(\Oo)$; this set is finite and
its cardinality $h(\Oo)=|\Cl(\Oo)|$ is the (left ideal) {\em class number} of $\Oo$. The \emph{discriminant} of an order $\Oo$ is by definition $\disc(\Oo)=\det(\tr(\gamma_{i}\ov\gamma_{j})_{i,j\leq 4})$  for  an $\mathcal{O}_F$-basis $\{\gamma_{1},\dots,\gamma_{4}\}$ of $\Oo$. The reduced norm $\nr(\mfI)$ of a lattice $\mfI$ is  the fractional $\mathcal{O}_F$-ideal generated by all elements $\nr(\gamma)$ with $\gamma \in \mfI$. The dual of a lattice $\mfI$ is the lattice
$\mfI^*=\{\gamma\in\B(F),\ \tr(\gamma\mfI)\subset \mathcal{O}_F\}.$ One defines the {\em reduced discriminant} $\rdisc(\Oo)$ of $\Oo$ to be  the ideal 
\begin{equation}\label{level}
  \rdisc(\Oo) := \nr((\Oo^{\ast})^{-1}) = \nr(\Oo^{\ast})^{-1}, 
\end{equation}  
so the reduced discriminant is the level of the associated norm  form.  
If $\Oo$ is a maximal order, its reduced discriminant equals the reduced discriminant $D_{\rmB}$ of $\B$  \cite[II.4.7]{Vig}.  In general, one has the following important relation between discriminant and reduced discriminant \cite[I.4.7]{Vig}
\begin{equation}\label{disc}\disc(\Oo)=\rdisc(\Oo)^2.\end{equation}

 An Eichler order is by definition the intersection of two maximal orders. To an Eichler order there is associated an $\mathcal{O}_F$-ideal $\mathfrak{N}$ coprime to $D_{\rmB}$ such that for every $\mathfrak{p}$ coprime with ${D}_{\B}$, $\Oo_{\mathfrak{p}}$ is conjugate to the order 
 \begin{equation}\label{order} \left(\begin{matrix} \mathcal{O}_{F, \mathfrak{p}} & \mathcal{O}_{F, \mathfrak{p}} \\ {\tt N}_\mathfrak{p} \mathcal{O}_{F, \mathfrak{p}} & \mathcal{O}_{F, \mathfrak{p}}\end{matrix}\right)\subset {\mathrm M}_{2}(F_{\mathfrak{p}})
 \end{equation}
   where $\tt N$ is a finite idele corresponding to the ideal $\mathfrak{N}$.   We may and will assume that the choice of $\Oo_{max}$ is such that $\Oo_{\mathfrak{p}}$ corresponds precisely to \eqref{order}.    Eichler orders associated with the same $\mathfrak{N}$ are locally conjugate (and conversely). Moreover,  the left order $\Oo_l(I)$ of a right $\Oo$-ideal $I$ is everywhere locally conjugate to $\Oo$, hence the norm forms of $\Oo_l(I)$ and $\Oo$ are in the same genus.
 
For an Eichler order  $\Oo$, the discriminant and the reduced discriminant have the following explicit expressions \cite[p.\ 85]{Vig} 
\begin{equation*}
\rdisc(\Oo)= \ D_{\B}\mathfrak{N},\quad \disc(\Oo)= (D_{\B}\mathfrak{N})^2,
\end{equation*}
and 
the class number equals \cite[p.\ 143]{Vig}
\begin{displaymath}
h(\Oo)\asymp  {\N}(D_{\B}\mathfrak{N})\prod_{\mathfrak{p}|D_{\B}}\left(1-\frac{1}{{\N}(\mathfrak{p})}\right)\prod_{\mathfrak{p}|\mathfrak{N}}\left(1-\frac{1}{\N(\mathfrak{p})}\right)^{-1}={\N(\disc(\Oo))}^{1/2+o(1)}.
\end{displaymath}

\subsection{Quaternary quadratic spaces and quaternion algebras}\label{sec61} 
Let $(V,Q)$ be a non-degenerate quaternary quadratic space over $F$ and  let $E$ be the quadratic $F$-algebra 
$$E=
\begin{cases}F\times F&\hbox{ if $\disc(Q)\in (F^\times)^2$ (i.e.\ is a square class),}\\
F(\sqrt{\disc(Q)})&\hbox{ if $\disc(Q)\not\in (F^\times)^2$,}
\end{cases}$$
equipped with either the $F$-invariant involution
$\sigma(x,y)=(y,x)\hbox{ if $E=F\times F$}$ or the canonical $F$-invariant involution if $E$ is a field. In the split case we view $F$ as embedded diagonally into $F\times F$. 

There is a unique quaternion algebra $\B$ defined over $F$ such that $(V,Q)$ is similar to the quaternary quadratic space $(\B',\Nr)$ that we now  describe (see also \cite{ponomarev}). Let  $\B_E:=\B\otimes_FE$ (this is often referred to as the second Clifford algebra of the quadratic space $V$). Slightly abusing notations, we denote by ${\cdot}^*$ the extension to $\B_E$ of the canonical involution of $\B$, by $\Nr(z)=z z^*$ the associated norm form on $\B_E$ and by $\sigma=\mathrm{Id}_\B\otimes\sigma$ the extension of $\sigma$ from $E$ to $\B_E$.
Let $$\B':=\{z\in \B_E,\ {\sigma(z)}^*=z\}.$$
Then $\Nr$ is $F$-valued on $\B'$ and $(\B',\Nr)$ defines a non-degenerate quaternary quadratic space over $F$ such that $\disc(\B')$ is a square if $E=F\times F$ and $E=F(\sqrt{\disc(\B')})$ otherwise.

We now proceed to describe the orthogonal group $\SO(\B')$ along these lines: for any $w\in \B_E^\times$, the map
$$z\mapsto wz\sigma(w)^*$$
leaves $\B'$ invariant and defines a proper similitude with factor $\lambda(w)=\Nr_{E/F}(\Nr(w))$. In particular, if $w$ is such that $\Nr(w)\in F$ (i.e. $\Nr(w)=\sigma(\Nr(w))$), the map
$$\rho_{w}: z\in\B'\mapsto wz\sigma(w)^{-1}=\frac{1}{\Nr(w)}wz\sigma(w)^* $$
is a special orthogonal transformation of $(\B',\Nr)$; moreover the map $w\mapsto \rho_w$ induces an isomorphism of $F$-algebraic groups
$$\SO(\B')\simeq\rmZ_F\bash\{ w\in \B^\times_E,\ \Nr(w)=\sigma(\Nr(w))\}.$$
Here we view $\B^\times_E$ as an $F$-algebraic group (of dimension $8$) and  $\Nr:\Bt_E\mapsto \Bbb{G}_{m,E}$  
and $\sigma:\Bt_E\mapsto \Bt_E$ as algebraic maps. We also note that the stabilizer in $\SO(\B')$ of the vector $1\in\B'$ is precisely
$$\SO(\B')_1=\rmZ_F\bash\{ w\in \B^\times_F\}=\mathrm{PB}^\times_F=\SO({\B'}^{0})$$
where ${\B'}^{0}=\B^0_E\cap\B'$ is the orthogonal subspace to $1$.
In the split case 
 we have  $\Bt_E=\Bt\times\Bt$, $\B'$ is identified with $\B$ via the embedding $z\in\B\mapsto (z,z^*)\in\B\times\B$, and the identification
$$\SO(\B)\simeq\rmZ_F\bash\{ (w,w')\in \Bt_F\times\Bt_F,\ \Nr(w)=\Nr(w'))\}$$
($Z$ diagonally embedded in $\Bt\times\Bt$) is given via the map
$$(w,w')\mapsto (\rho_{w,w'}:z\mapsto wz{w'}^{-1}).$$
We denote by $\rmG<\tG$ the $F$-algebraic groups defined (at the level of their $F$-points) by
\begin{equation}
\label{Gdef}
\rmG=\rmZ_F\bash\{ w\in \B^\times_E,\ \Nr(w)=\sigma(\Nr(w))\}\simeq\SO(\B'),\quad  \tG=\rmZ_F\bash\B^\times_E.
\end{equation}

\subsection{Representations at the archimedean place}\label{archi_preliminaries}

For any integer $m\geq 0$ there is a unique irreducible (unitary) representation of $\SU\simeq\mathrm{Spin}_3(\Rr)$ of degree $d_m=m+1$, denoted $\pi_{m}$, and any irreducible representation of $\SU$ is isomorphic to some $\pi_m$, see e.g.\ \cite[Section 7]{Fa} for details.  The representation $\pi_m$  may be  realized concretely as the space of complex homogeneous polynomials of degree $m$ in two variables on which $\SU\subset\GL_2(\Cc)$ acts by linear change of variables. The Casimir element (say with respect to the inner product on the Lie algebra $(X,Y)=-\frac12\mathrm{Tr}(XY)$) $C_{\SU}\in {U}(\mathfrak{su}_2)_\Cc$ acts  on any realization by multiplication by the scalar $$\lambda_m=-{m}({m}+2).$$

Let $\TR \cong \SO_2(\Bbb{R}) \subset\SU$ be the stabilizer of (say) the north pole of $S^2$ under the natural projection $\SU\ra\SO_3(\Rr)$.  This is a maximal torus of $\SU$ isomorphic to $\SO_2(\Rr)\simeq \Rr/2\pi\Zz$ and whose image in $\SO_3(\Rr)$ is the group of matrices
$$\SO_2(\Rr)=\biggl\{\left(\begin{array}{ccc}\cos(2\theta) & -\sin(2\theta) & 0 \\\sin(2\theta) & \cos(2\theta) & 0 \\0 & 0 & 1\end{array}\right),\ \theta\in\Rr/2\pi\Zz\biggr\}.$$
Let $\theta\in\Rr/2\pi\Zz\ra \kappa(\theta)\in\TR$ be a parametrization and $e:\TR\ra\Cc^1$ be the character
$$e(\kappa(\theta))=\exp(\iota\theta).$$ 
If $V_m$ is any vector space realizing $\pi_m$  and $l\in\Zz$, let $V_m^l$ be the subspace of vectors ``of weight $l$", that is, the vectors satisfying
\begin{equation*}
\vphi\in V_m, \quad \kappa.\vphi=e(\kappa)^l\vphi \text{ for all }\kappa\in\TR.
\end{equation*}
Then $V_m^l$ is one dimensional if $|l|\leq m,\ l\equiv m\ ( 2)$ and zero otherwise. 
\begin{rem}The representation $\pi_m$ occurs in the right regular representation $L^2(\SU)$ with multiplicity $d_m$. When $m$ is even, $\pi_m$ descends via the natural projection $\SU\rightarrow\SO_3(\Rr)$ to an irreducible representation of $\SO_3(\Rr)$. The direct sum of the weight zero vectors of each such copy of $\pi_m$ therefore injects into $L^2(\SO_3(\Rr))^{\SO_2(\Rr)}=L^2(S^2)$ and the image is the space of harmonic homogenous polynomials of degree $m/2$ in $\Rr^3$ (i.e.\ polynomials $P$  such that $\Delta_{\Rr^3}P=0$). The action of the Casimir element on this space corresponds to that of a fixed multiple of the Laplace operator $\Delta_{S^2}.$
\end{rem}

Given a non-zero $\vphi\in V_m^l$   and $g\in\SU$, we write
$$p_{m, l}(g):=\frac{(\vphi,g.\vphi)_m}{(\vphi,\vphi)_m}$$
for the corresponding normalized matrix coefficient, where $(\ ,\ )_m$ is some $\SU$-invariant inner product  on $V_m$. By definition $g\mapsto |p_{m, l}(g)|$ is bi-$\TR$-invariant, and therefore depends only on $$t=t(g)=\peter{\rho_g(x_3),x_3}_{\Rr^3} \in [-1, 1],$$
the inner product of the north pole $x_3$ on $S^2$ with its image by the corresponding rotation. The following  decay estimate holds as $g\in\SU$ gets ``away" from $\TR$ (i.e.\ $t$ gets away from $\pm 1$). We will use it for the bound in the eigenvalue aspect in  Section \ref{eig-aspect}. 

\begin{lemma}\label{decaylemma}
For $-m\leq l\leq m$ one has
\begin{equation} \label{decayeqn} |p_{m,l}(g)|\ll \min\left(1, \Bigl(\frac{m+1}{|l|+1}\Bigr)^{-1/2}(1-t^2)^{-1/4}\right).
\end{equation}
\end{lemma}
\proof By symmetry we may  assume  $0\leq l\leq m$.  Being a matrix coefficient, it is clear that $|p_{m,l}(g)|\leq 1$. One has \cite[\S 6.3.1 \& 6.3.7]{Vil}
$$|p_{m,l}(g)|=|P^m_{l,l}(t)|=\Bigl(\frac{1+t}2\Bigr)^{l}|P^{(0,2l)}_{m-l}(t)|$$
where (for $\alpha,\beta\geq 0$ integers)
$$P_{n}^{(\alpha,\beta)}(t)=\frac{(-1)^{n}}{2^{n}n!}(1-t)^{-\alpha}(1+t)^{-\beta}\frac{d^{n}}{dt^{n}}[(1-t)^{\alpha}(1+t)^{\beta}(1-t^2)^{n}]$$
 is the Jacobi polynomial. Let us recall that $P_{n}^{(\alpha,\beta)}$ has degree $n$ and that $\{P_{n}^{(\alpha,\beta)}\mid n\geq 0\}$ is orthogonal with respect to the inner product 
$$\peter{P,Q}_{(\alpha,\beta)}:=\int_{-1}^1P(t)\overline{Q(t)}(1-t)^{\alpha}(1+t)^{\beta}dt,$$
and that \cite[6.10.1(7)]{Vil}
$$\peter{P_{n}^{(\alpha,\beta)},P_{n}^{(\alpha,\beta)}}_{(\alpha,\beta)}=\frac{2^{\alpha+\beta+1}}{2n+\alpha+\beta+1}\frac{(n+\alpha)!(n+\beta)!}{(n+\alpha+\beta)!n!}.$$
In particular, 
$$\peter{P_{m-l}^{(0,2l)},P_{m-l}^{(0,2l)}}_{(0,2l)}=\frac{2^{2l+1}}{2m+1}\frac{(m-l)!(m+l)!}{(m+l)!(m-l)!}=\frac{2^{2l+1}}{2m+1}. $$
By \cite[Thm.\ 1]{EMN} we have 
$$(1-t)^{\alpha/2}(1+t)^{\beta/2}P_{n}^{(\alpha,\beta)}(t)\ll \peter{P_{n}^{(\alpha,\beta)},P_{n}^{(\alpha,\beta)}}_{(\alpha,\beta)}^{1/2}\frac{(\alpha+\beta+1)^{1/2}}{(1-t^2)^{1/4}}$$
for $t\in[-1,1]$ and hence
$$(1+t)^l2^{-l}P_{m-l}^{(0,2l)}(t)\ll \frac{(l+1)^{1/2}}{(1-t^2)^{1/4}(2m+1)^{1/2}}.$$

\qed
\begin{rem}\label{remEMN} The above bound exhibits significant decay as $t$ gets away from $\pm 1$ uniformly for $|l|\leq (1+m)^{1-\delta}$ for any fixed $\delta>0$. 
It is plausible that this holds also for very large values of $l \asymp m$: for instance in the extreme case  $l=m$ one has:
$$|p_{m,m}(g)|=\Bigl(\frac{1+t}2\Bigr)^m.$$
More generally it is conjectured in  \cite{EMN} that in \refs{decayeqn}, the term $\frac{m+1}{|l|+1}$ can be replaced by $\frac{m+1}{\sqrt{|l|+1}}$.
\end{rem}

We extend these notation to irreducible representations of products: for 
$\bfm=(m_\sigma)_\sigma$ a $d$-tuple of non-negative integers, we set $$\pi_\bfm=\bigotimes_{i=1}^d \pi_{m_i}\in\mathrm{Irr}(\SU^d)\hbox{ and denote by
 }d_\bfm=\prod_\sigma ({m_\sigma}+1)=:|\bfm|$$
its dimension (cf.\ footnote 1). 
Given a realization $V_{\bfm}=\bigotimes_{\sigma}V_{m_\sigma}$   of $\pi_\bfm$ and a $d$-tuple $\bfl=(l_\sigma)_\sigma$   of integers, we denote by $V_\bfm^\bfl=\bigotimes_{\sigma}V_{m_\sigma}^{l_\sigma}$ the tensor product of weight $l_\sigma$ vectors with respect to the product $\TR^d$, and for $g_\infty=(g_\sigma)_\sigma\in\SU^d$ we denote by
\begin{equation}\label{matrix}
 p_{\bfm, \bfl}(g_\infty):=\prod_\sigma p_{m_\sigma, l_{\sigma}}(g_\sigma)
\end{equation} 
  the corresponding normalized matrix coefficient.

\subsection{Adelic interpretation of ellipsoids}\label{adelicinterp}
 As in \cite[\S 4]{BM1}, we now define $X^{(2)}(\Oo)$ as an adelic quotient. We refer to \cite[Chap. III \& V]{Vig} for more details.
 
\subsubsection{Archimedean place}\label{archisubsub} Fix once and for all $\{x_1,x_2,x_3\}=(\{x_{1,\sigma},x_{2,\sigma},x_{3,\sigma}\})_\sigma$ an orthonormal basis of $\B^0(F_\infty)=\prod_{\sigma} \B^0(F_\sigma)\simeq (\Rr^3)^d$ for the quadratic form $\oplus_\sigma \nr$; this induces identifications $$\B^1(F_\infty)\simeq\SU^d,\ \PB(F_\infty)\simeq\SO_3(\Rr)^d.$$
 Let $$x^0:=x_3=(x_{3,\sigma})_\sigma\in\B^{0,1}(F_\infty)\simeq (S^2)^d\hbox{ and }\rmK_\infty\simeq\SO_2(\Rr)^d$$ be its stabilizer under the conjugacy action of $\B^1(F_\infty)$ on $\B^0(F_\infty)$; this yields an identification 
 \begin{equation}\label{archiident}
 \rmZ(F_\infty)\bash\Bt(F_\infty)/\rmK_\infty\simeq (S^2)^d. 
 \end{equation}
 \subsubsection{Finite places}\label{finitesubsub}
Given an Eichler order $\Oo$ (the intersection of two maximal orders), we denote by $\widehat{\Oo}$ the closure of $\Oo$ in $\B(\Aa_f)$. Correspondingly let $K_f=\whOt\subset \Bt(\Af)$ be the subgroup of units  of the ring $\whO$ ($K_f$ is an open compact subgroup of $\Bt(\Af)$) and $\whO^0\subset\B^0(\Af)$ the elements of $\whO$ with trace 0. 

We define the manifold in Theorem \ref{mainthm} as the adelic quotient   
\begin{equation}\label{defX2adelic}
  X^{(2)}(\Oo):= Z(\Aa)\Bt(F)\bash \Bt(\Aa)/\rmK_\infty \whOt.
\end{equation} 
Because of \eqref{archiident},  $X^{(2)}(\Oo)$ is therefore identified (cf. \refs{defX2}) with a finite disjoint union of quotients of $d$-fold product of spheres  indexed by the finite double coset $Z(\Af)\Bt(F)\bash \Bt(\Af)/\whOt$. Alternatively one could also think of $X^{(2)}(\Oo)$ as a collection of quotients of ellipsoids associated to different quadratic forms (covering all classes in the genus of $(\Oo^0,\nr)$). 

Moreover, if $\Oo$ is assumed to be an Eichler order, every  $\Oo$-ideal in $\B(F)$ is locally principal  and the double cosets $Z(\Af)\Bt(F)\bash \Bt(\Af)/\whOt$ parametrize precisely $\mathrm{Cl}(\Oo)$.

\subsubsection{Measures}\label{measures} We set $$[\Bt]:=Z(\Aa)\Bt(F)\bash \Bt(\Aa);$$ this quotient
(equipped with the quotient of Haar measures, see \cite[Cor. 2.3, Chap.\ V]{Vig}) has finite volume which is equal to $2$ when the chosen measures are the Tamagawa measures. The volume of $\whOt\subset\Bt(\Af)$ for the Tamagawa measure is of size $\asymp_F \N(\disc(\Oo))^{-1/2+o_F(1)}$; therefore, if we normalize Haar measures on $\Bt$ so as to match the normalization given in  the introduction at the infinite places, and for the finite places so that $\whOt$ has mass $1$, and denote the resulting measure on $[\Bt]$ by $dg$, we obtain
$$\int_{[\Bt]}dg=\vol(X^{(2)}(\Oo))=V_2$$
which is precisely \refs{defV2}. The corresponding inner product will be denoted by
$$\peter{\vphi,\vphi}=\int_{[\Bt]}|\vphi(g)|^2dg.$$

\subsubsection{Automorphic forms} Our original problem is equivalent to bounding a certain $L^2$-normalized {\em automorphic function} $\vphi$ on the adelic quotient $Z(\Aa)\Bt(F)\bash \Bt(\Aa)$, which is
\begin{enumerate}
\item right $\rmK_\infty.\rmK_f$-invariant,
\item an  eigenfunction of the Casimir operators $C=(C_\sigma)_\sigma$ associated to the group $\B^1(F_\infty)\simeq \mathrm{Spin}_3(\Rr)^d$ with eigenvalues 
\begin{equation}\label{laplaceeigen}
\lambda_\sigma=-{m_\sigma}({m_\sigma}+2)
\end{equation}
for a $d$-tuple $\bfm$ of even integers as in the preceding subsection, 
\item an  eigenfunction of a certain {\em Hecke algebra} $\mcH(\Oo)$, which is a commutative algebra of normal operators commuting with the Casimir operators $C$.
\end{enumerate}
We will recall the definition of the Hecke algebra $\mcH(\Oo)$ in section \S \ref{pretrace}; in particular we do not consider Hecke operators at ramified primes, so that the operators of $\mcH(\Oo)$ are indeed normal.  Our assumptions imply that the $\BtA$-translates of $\vphi$ generate an automorphic representation $\pi=\otimes_v\pi_v$ of $\BtA$ with trivial central character such that $$\pi_\infty=\otimes_{\sigma}\pi_\sigma\simeq\pi_\bfm$$ for $\bfm=(m_\sigma)_\sigma$ given in \refs{laplaceeigen} and $\pi_\bfm$ defined in \S \ref{archi_preliminaries} (under the identifications of \S \ref{archisubsub}), and such that $\pi_v$ is an unramified representation for every finite place $v\nmid\disc(\Oo)$.

 If $\pi$ is finite-dimensional, then it is one-dimensional and $\vphi$ is proportional to the function
$$g\in\BtA \ra \chi(\Nr(g))$$ for   some  (quadratic) character $\chi$ on $F^\times\bash\At$. In this case  $\vphi$  is constant on the various components of $X^{(2)}(\Oo)$ with value equal to $\pm V_2^{-1/2}$.  In particular its Laplace eigenvalues are $(0,\cdots,0)$, and the bounds of Theorem 1 are a fortiori satisfied. Therefore we can restrict ourselves to infinite dimensional representations for the rest of the paper.

 \subsection{A slight generalization}
With no extra effort we can consider a slightly more general setting: let $\chi:F^\times\bash\At\ra \Cc^1$ be a unitary Hecke character, and $L^2(Z(\Aa)\Bt(F)\bash\BtA, \chi)$
the space functions on $\Bt(F)\bash\BtA$ satisfying $$\vphi(\lambda \gamma g)=\chi(\lambda)\vphi(g),\ \lambda\in Z(\Aa),\ \gamma\in\Bt(F),\ g\in\BtA$$ and such that $\peter{\vphi,\vphi}=\int_{[\Bt]}|\vphi(g)|^2dg$ is finite.
 
Let 
 $\pi\simeq\otimes_v\pi_v\hookrightarrow L^2(Z(\Aa)\Bt(F)\bash\BtA,\chi)$ be an infinite dimensional irreducible automorphic representation of $\BtA$ with central character $\chi$. 
 Under the identification $\B^1(F_\infty)\simeq \SU^d$, the representation $\pi_\infty:=\otimes_{\sigma\mid\infty}\pi_\sigma$ corresponds  to $\pi_\bfm$ for some $\bfm\in\Bbb{N}_{\geq 0}^d$.  Let 
 \begin{equation*}  
\lambda=(\lambda_\sigma)_\sigma,\ \lambda_\sigma=-m_\sigma(m_\sigma+2)
\end{equation*}
 be the eigenvalues of the Casimir operators $(C_{\B^1(F_\sigma)})_\sigma$ and let as in \eqref{deflambda}
 $$|\lambda|=\prod_\sigma(1+|\lambda_\sigma|).$$
 
 Let $\Oo\subset\B$ be an Eichler order and let $\vphi\in\pi$ be a non-zero smooth, $\whOt$-invariant function of some fixed weight   $\bfl=(l_\sigma)_\sigma\in\Zz^d$ with respect to action of  the maximal torus $\rmK_\infty\simeq \SO_2(\Rr)^d$ (in particular $\pi_v$ is an unramified principal series representation at every place not dividing $\disc(\Oo)$ and $\chi$ is unramified at these places as well).  

Under these conditions, we prove the following slightly more general version of  Theorem \ref{mainthm}:
 \begin{thm}\label{mainthm2} Let $\Oo$ be an Eichler order in a totally definite quaternion algebra $\B$ over $F$, and let $\vphi \in \pi \subseteq L^2(Z(\Aa)\Bt(F)\bash\BtA,\chi)$ as described in this subsection. Then one has
 \begin{equation*}
\|\vphi\|_\infty \ll_\bfl |\lambda|^{\frac{1}{4}}  (V_2|\lambda|^{\frac{1}{2}})^{-\frac{1}{20}}\|\vphi\|_2.
\end{equation*}
We obtain the following individual bounds in the $\lambda$ and in the volume aspect
$$\|\vphi\|_\infty \ll_\bfl |\lambda|^{\frac{1}{4}} V_2^{-\frac{1}{6}+\varepsilon}, \quad \ \|\vphi\|_\infty \ll_\bfl |\lambda|^{\frac{1}{4} - \frac{1}{27}} \|\vphi\|_2.$$
All implied constants depend at most on $F$, $\varepsilon$ and $\bfl$. 
\end{thm}
The next three sections are devoted to the proof of this theorem.
 \begin{rem} Here we have assumed that the weight $\bfl=(l_\sigma)_\sigma$ of $\vphi$ is fixed. This is merely to simplify exposition (and also because the lowest weight case is arguably the most interesting one). The proof of Theorem \ref{mainthm2} together with the bound \refs{decayeqn} yields immediately a non-trivial bound for $\|\vphi\|_\infty$ uniformly across all weights $\bfl$ satisfying 
 $$|\bfl|=\prod_\sigma(1+|l_\sigma|)\leq |\bfm|^{1-\delta}$$
 for any fixed $\delta>0$. We   expect that such a non-trivial bound holds for all $\bfl$, and this would follow from good enough
 bounds for Jacobi polynomials (cf. Remark \ref{remEMN}). 
 \end{rem}

\section{Reduction of definite quadratic forms}\label{reduction}

Let
\begin{displaymath}
  Q(\textbf{x}) =  \frac{1}{2} \sum_{1 \leq i, j \leq n} a_{ij} x_i x_j, \quad a_{ij} = a_{ji} \in \mathcal{O}_F, \quad a_{jj} \in 2 \mathcal{O}_F,
\end{displaymath}
be an $F$-integral   quadratic form in $n$ variables. Let $A = (a_{ij})_{1 \leq i, j \leq n} \in {\rm Mat}_{n\times n}(\mathcal{O}_F)$ be the symmetric $n\times n$-matrix associated to $Q$. The {\em determinant} of $Q$ is 
$$\Delta=\det A. $$Note that the determinants of two equivalent forms over $\mathcal{O}_F$ may differ by the square of a unit. The quadratic form defines a bilinear form
\begin{equation}\label{inner}
 \langle\textbf{x}, \textbf{y}\rangle := \frac{1}{2} \textbf{x}^t A \textbf{y} = \frac{1}{2}(Q(\textbf{x} + \textbf{y}) - Q(\textbf{x}) - Q(\textbf{y})). 
\end{equation}
Let $\mathcal{O}_F^{\sharp} := \{\textbf{y} \in F^n \mid 2\langle\textbf{y}, \mathcal{O}_F\rangle \subseteq \mathcal{O}_F\}$. The {\em level} of $Q$ is the integral  ideal $\mathfrak{n} := (Q(\mathcal{O}_F^{\sharp})\mathcal{O}_F)^{-1}$. In particular, if $N \in \mathfrak{n}$, then $NA^{-1} \in {\rm Mat}_{n\times n}(\mathcal{O}_F)$ is an integral matrix. (Indeed, if $N \in \mathfrak{n}$, then by definition $N \textbf{x}^t A \textbf{x} \in 2 \mathcal{O}_F$ for all $\textbf{x} \in F^n$ such that $\textbf{x}^t A \textbf{y} \in \mathcal{O}_F$ for all $y \in \mathcal{O}_F^n$, hence $\textbf{z}^t (NA^{-1})\textbf{z} \in 2 \mathcal{O}_F$ for all $\textbf{z} \in \mathcal{O}_F^n$ which implies that $NA^{-1}$ is integral.)

For any real embedding $\sigma:F\hookrightarrow \Rr$ denote by $Q^{\sigma}$ the conjugated form. We assume that $Q$ is totally positive definite, that is,  $Q^{\sigma}$ is positive definite for all $\sigma$.

Minkowski developed a reduction theory for rational positive definite quadratic forms (see e.g.\ \cite[Chapter 12]{Cas}) that has been extended to arbitrary number fields by Humbert \cite{Hum}. We summarize some basic facts. Every quadratic form is equivalent (over $\mathcal{O}_F$) to some form of the shape
 \begin{equation}\label{quasidiag}
 \begin{split}
 Q(\textbf{x}) & = \frac{1}{2}\sum_{1 \leq i, j \leq n} a_{ij} x_ix_j \\
 &= h_1 (x_1 + c_{12} x_2 + \ldots + c_{1n} x_n)^2 +  h_2(x_2 + c_{23} x_3 + \ldots + c_{2n}x_n)^2 + \ldots + h_nx_n^2
 \end{split}
\end{equation}
with $c_j, h_j \in F$ where 
\begin{equation}\label{sizea}
   a_{ij}^{\sigma} \ll a_{jj}^{\sigma} \asymp h_j^{\sigma}
\end{equation}   
    for all $1 \leq i, j \leq n$ and all embeddings $\sigma$, and 
\begin{equation}\label{sizeh}
 1\ll  h^{\sigma_1}_1 \asymp \ldots \asymp h^{\sigma_{d}}_1 \ll  h^{\sigma_1}_2 \asymp \ldots \asymp h^{\sigma_{d}}_2 \ll  \ldots \ll  h^{\sigma_1}_n \asymp \ldots \asymp h^{\sigma_{d}}_n.
\end{equation}
Here and henceforth all implied constants depend only on $n$ and $F$. Clearly, 
\begin{equation}\label{allh}
  h_1^{\sigma} \cdot \ldots \cdot h_n^{\sigma} = \Delta^{\sigma}
\end{equation}
where (by slight abuse of notation) $\Delta$ is the determinant of the form \eqref{quasidiag}. This  determinant (which may differ from the determinant of the original form $Q$) has the advantage that its conjugates are of comparable size. From now on we will always refer to this balanced determinant when we use the symbol $\Delta$. Of course, this convention is not necessary when we use the norm of $\Delta$.

We denote the eigenvalues of the  matrix $A^{\sigma} = (a^{\sigma}_{ij})$ by $0 < \lambda^{\sigma}_1 \leq \lambda^{\sigma}_2 \ldots \leq \lambda^{\sigma}_n$. By \eqref{sizea} -- \eqref{allh}, the determinant of any $(n-1) \times (n-1)$ submatrix of $A^{\sigma}$ is   $O(\Delta^{\sigma})$, hence by Cramer's rule the   eigenvalues of $(A^{\sigma})^{-1}$ are $O(1)$, and therefore 
 \begin{equation}\label{boundeig}
  1 \ll \lambda^{\sigma}_1 \leq \lambda^{\sigma}_n \ll \Delta^{\sigma}.  
\end{equation}
Let $\tilde{Q}$ be the quadratic form in $n-1$ variables that is derived from $Q$ by setting $x_n = 0$.  Let $\tilde{A}$  be the corresponding  $(n-1)\times (n-1)$-submatrix of $A$, and denote by $\tilde{\Delta}$ its determinant. The $(n, n)$th-entry  of $A^{-1}$ is by Cramer's rule $\tilde{\Delta}/\Delta$ (up to sign); hence $N\tilde{\Delta}/\Delta$ is integral for all $N \in \mathfrak{n}$. Therefore the ideal $\mathfrak{n}\, (\tilde{\Delta})/(\Delta)$ is integral, and we obtain
\begin{equation}\label{dtilde}
 \N (h_1)  \times \ldots \times \N (h_{n-1})   =  \N (\tilde{\Delta})  \geq \frac{\N(\Delta)}{\N(\mathfrak{n})}.
\end{equation}  

\section{Representation numbers of quadratic forms}\label{sec:quadratic}

In this section we establish several lemmata to bound certain averages of  representation numbers of $\mathcal{O}_F$-integers  by some totally definite quadratic form $Q$.   To perform the counting we will frequently use the following  consequence of Dirichlet's unit theorem: let $A_1, \ldots, A_d > 0$ be any positive real numbers and write $A = A_1 \cdot \ldots \cdot A_d$.  Then
\begin{equation}\label{units}
\#\{u \in U : |u^{\sigma_j}| \leq A_j\} \ll_F \log(2+A)^{d-1}. 
\end{equation}
As a consequence we find
 \begin{equation}\label{lattice}
  \#\{   x \in \mathcal{O}_F :  0 <  |x^{\sigma_j}|  \leq A_j\} \ll_F A.
  \end{equation}
Indeed,  \eqref{units} implies that for each principal ideal $(x)$ of norm $\N x \leq A$  there are $O(\log(2 + A/\N x)^{d-1})$ generators satisfying the size constraints in \eqref{lattice}, hence the left hand side of \eqref{lattice} is at most
\begin{displaymath}
  \ll \sum_{\N(\mathfrak{a}) \leq A} \log\left(2 + \frac{A}{\N(\mathfrak{a})}\right)^{d-1} \ll A. 
\end{displaymath}
We remark that the estimate \eqref{lattice} is a trivial lattice point count if all $A_j \gg 1$. It is a little less trivial if some $A_j$ are very large and others are very small. 

We use the notation $r_Q(\ell)$ to  denote  the number of integral representations of $\ell$ by $Q$.  From now we consider quadratic forms in 2, 3 and 4 variables. The following lemma derives uniform estimates for representation numbers of quadratic forms, averaged over thick and not so thick sequences (the thinner sequences are needed because of the special form of our amplifier). The lemma will be used in Section \ref{volume} when we derive the volume bound.

\begin{lemma}\label{lemma2} Let $Q$ be a totally positive-definite integral quaternary quadratic form of determinant $\Delta$ and level $\mathfrak{n}$. Let $y, y_1, y_2 >1$.  
Then
\begin{equation}\label{e3}
  \sum_{\substack{\ell \in \mathcal{O}_F \\ 0 \leq   \ell^{\sigma} \leq y^{1/d} }}   r_Q(\ell ) \ll \frac{y^2}{\N(\Delta)^{1/2}} + \frac{y^{3/2}}{(\N(\Delta)/\N(\mathfrak{n}))^{1/2}} + y.
\end{equation}
If in addition $h_1 \asymp 1$ in \eqref{quasidiag}, then 
\begin{equation}\label{e2}
  \sum_{\substack{ \ell_1 \in \mathcal{O}_F\\  0 \leq \ell_1^{\sigma} \leq y_1^{1/d}}} \sum_{\substack{\ell_2 \in \mathcal{O}_F \\ 0 \leq \ell_2^{\sigma} \leq y_2^{1/d}}}  r_Q(\ell_1 \ell_2^2) \ll y_1 \left(\frac{(y_1y_2^2)^{3/2}}{\N(\Delta)^{1/2}} + \frac{y_1y_2^2}{(\N(\Delta)/\N(\mathfrak{n}))^{1/2}} + (y_1y_2^2)^{1/2}\right)(y_1y_2\N(\Delta))^{\varepsilon}, 
\end{equation}
\begin{equation}\label{e1}
  \sum_{\substack{\ell\in \mathcal{O}_F\\ 0 \leq \ell^{\sigma} \leq y^{1/d}}} r_Q(\ell^2) \ll \left(\frac{y^3}{\N(\Delta)^{1/2}} + \frac{y^{2}}{(\N(\Delta)/\N(\mathfrak{n}))^{1/2}} + y\right)(y\, \N(\Delta))^{\varepsilon} 
\end{equation}
for any $\varepsilon>0$, the implied constants depending on $\varepsilon$ alone. Here and in the following a summation condition of the type  $0 \leq   \ell^{\sigma} \leq y^{1/d}$ is understood to hold for all embeddings $\sigma$. 
\end{lemma}

The first and last term on the right hand side of \eqref{e3} and \eqref{e1} are certainly   optimal. Maybe the middle term can be improved slightly. The bound \eqref{e2} is not best possible in general, but sufficient for our purposes. \\ 

\textbf{Proof.} All of these bounds are proved in a similar way. We start with \eqref{e3}. We use the representation \eqref{quasidiag} together with the bounds \eqref{sizeh}. Let $\tilde{x}_j := x_j + \sum_{i > j} c_{ji} x_i$.  By \eqref{lattice} we have  $\ll (y/\N (h_j))^{1/2}$ non-zero choices for $\tilde{x}_j$. Hence we have $\ll (y/ \N (h_j))^{1/2} + 1$ choices in total for $\tilde{x}_j$, and hence  $\ll (y/ \N (h_j))^{1/2} + 1$ choices for $x_j$ getting a bound
\begin{gather*}
 \ll \left(\frac{y}{\N (h_4)} + 1\right)^{1/2}\left(\frac{y}{\N (h_3)} + 1\right)^{1/2}\left(\frac{y}{\N (h_2)} + 1\right)^{1/2}\left(\frac{y}{\N (h_1)} + 1\right)^{1/2}\\ \ll \frac{y^2}{\N \Delta^{1/2}} + \frac{y^{3/2}}{(\N \Delta/\N(\mathfrak{n}))^{1/2}} + y
\end{gather*}
by \eqref{sizeh}, \eqref{allh} and  \eqref{dtilde}. \\

In order to prove \eqref{e1}, we choose as before $x_4, x_3, x_2$ in 
\begin{displaymath}
 \ll \left(\frac{y^2}{\N (h_4)} + 1\right)^{1/2}\left(\frac{y^2}{\N (h_3)} + 1\right)^{1/2}\left(\frac{y^2}{\N (h_2)} + 1\right)^{1/2}  \ll \frac{y^3}{\N \Delta^{1/2}} + \frac{y^2}{(\N(\Delta)/\N(\mathfrak{n}))^{1/2}} + y
 \end{displaymath}
 ways. Here we used again \eqref{sizeh}, \eqref{allh} and  \eqref{dtilde}. Note that $Q(\textbf{x}) = \ell^2$ implies  $x_j^{\sigma} \ll y^{1/d}$ for all $\sigma$, since $\lambda_1^{\sigma} \gg 1$ by \eqref{boundeig}.  Once we have fixed $x_2, x_3, x_4$ we are left with counting pairs $(x_1, \ell)$ satisfying
\begin{displaymath}
  2a_{11}\ell^2 - (a_{11}x_1 + \xi)^2 = D
\end{displaymath} 
 where 
\begin{displaymath}
  \xi =  \sum_{j=2}^4 a_{1j} x_j, \quad D =   - \sum_{i, j = 2}^4 (a_{1i}a_{1j} - a_{11}a_{ij}) x_ix_j . 
\end{displaymath}
Note that $\xi^{\sigma} \ll  y^{1/d}$ by \eqref{sizea} and our assumption $h_1 \asymp 1$,  and $D^{\sigma} \ll (\Delta^{\sigma}) y^{2/d}$ by the same argument. It follows that   $x_1^{\sigma} \ll  (\Delta^{\sigma})^{1/2} y^{1/d}$ (although precise exponent are insignificant here -- we only need polynomial dependence). 
 
Let us first assume that $2a_{11} = b^2$, say,  is a square in $F$ and hence in $\mathcal{O}_F$.   If $D\not=0$, then by a standard divisor argument there are   $\ll (\N D)^{\varepsilon}$ pairs of principal ideals $(b\ell - a_{11}x_1 -\xi)$, $(b\ell + a_{11}x_1 +\xi)$ whose product equals $(D)$, and by \eqref{units} each of these has   $\ll (y\, \N(\Delta))^{\varepsilon}$ generators $g\in \mathcal{O}_F$ satisfying $g^{\sigma} \ll (\Delta^{\sigma})^2 y^{2/d}$. This in turn gives  $\ll (y\, \N(\Delta))^{\varepsilon}$ choices for $x_1$. 
 
If $D=0$, we choose $\ell$ freely in $O(y)$ ways (by \eqref{lattice}), and then there are at most two choices for $x_1$. We determine how often the case $D=0$ happens. The quantity $D$ is a \emph{definite} ternary  quadratic form in $x_2, x_3, x_4$ whose determinants of its upper left $k \times k$ submatrices ($1 \leq k \leq 3$) are precisely the determinants of the $(k+1) \times (k+1)$ upper left submatrices of $A$ (up to sign). In particular we see that $D = 0$ if and only if $x_2 = x_3 = x_4 = 0$.

Let us now assume that $2a_{11}$ is not a square in $F$. Then $D \not=0$,  and we need to solve a Pell-type equation. There are   $\ll (\N D)^{\varepsilon}$ ideals $(\sqrt{2a_{11}}\ell  - a_{11} x_1 - \xi )$ in the totally real field $E = F(\sqrt{2a_{11}})$ of relative norm $D$, and again by \eqref{units} each of these yield  $\ll (y \, \N(\Delta))^{\varepsilon}$ solutions for $x_1$.  This establishes \eqref{e1}.\\


Finally we prove \eqref{e2}. Again we fix $x_4, x_3, x_2$ as above, and we fix $\ell_1$. This gives a total count of
 \begin{displaymath}
   y_1 \left(\frac{(y_1y_2^2)^{3/2}}{\N(\Delta)^{1/2}} + \frac{y_1y_2^2}{(\N(\Delta)/\N(\mathfrak{n}))^{1/2}} + (y_1y_2^2)^{1/2}\right),
 \end{displaymath}
 and we are left with counting pairs $(x_1, \ell_2)$ satisfying $2a_{11}\ell_1\ell_2^2 - (a_{11}x_1 + \xi)^2 = D $ with $\xi$ and $D$ as above. Now we argue exactly as in the previous case. \qed\\ 
 

For a polynomial 
\begin{displaymath}
  P(x_1, \ldots x_k) = \sum_{\textbf{n} \in \Bbb{N}_{\geq 0}^k} \alpha_{\textbf{n}}\textbf{x}^{\textbf{n}} \in \mathcal{O}_F[x_1, \ldots, x_k]
\end{displaymath}  
let $H(P ) := \sum_{\sigma}\sum_{ \textbf{n}} |\alpha_{\textbf{n}}^{\sigma}|$ denote its ``height". For $\ell \in \mathcal{O}_F$ we write $|\ell|_{\infty} := \max_{\sigma}|\ell^{\sigma}|$. The next two lemmas have preparatory character and are used as in input for the important Lemma \ref{VdKlemma} below. However, in particular Lemma \ref{lemma3} may be of independent interest, as it bounds \emph{uniformly} representation numbers of positive definite ternary and quaternary forms over (fixed) totally real number fields. 

\begin{lemma}\label{lemma3} a) Let $P(x, y) \in \mathcal{O}_F[x, y]$ be a quadratic polynomial and assume that its quadratic homogeneous part is a totally positive definite quadratic form. 
 Let $\ell \in \mathcal{O}_F$.  Then there are   $ \ll (H(P )(1+|\ell|_{\infty} )^{\varepsilon})$ solutions to $P(x, y)   = \ell$.

b) Let $Q$ be a totally positive definite integral ternary quadratic form over $F$ of discriminant $\Delta$ and let $\ell \in \mathcal{O}_F \setminus \{0\}$. 
Then $r_Q(\ell) \ll \N(\ell)^{1/2}  \, (| \ell|_{\infty}\N(\Delta))^{\varepsilon}$.


c) Let $Q$ be a totally positive definite integral quaternary quadratic form over $F$ of discriminant $\Delta$ and  let $\ell \in \mathcal{O}_F \setminus \{0\}$. 
Then $r_Q(\ell) \ll \N(\ell) \, (|\ell|_{\infty}\N(\Delta))^{\varepsilon}$. 

Here all implied constants depend on $\varepsilon$ at most. 
\end{lemma}

\begin{rem} The proof gives slightly stronger bounds for parts b) and c); for instance in the situation of part b) we obtain $r_Q(\ell) \ll (\N(\ell)^{1/2}\N(\Delta)^{-1/6}+1)( |\ell|_{\infty}\N \Delta)^{\varepsilon}$, but we do not need these refinements. We use part a) of the lemma in parts b) and c). Although they could also be proved without recourse to the first part, the first part will be needed in the proof of Lemma \ref{VdKlemma}, and so we take the opportunity to state it and prove it here. 
\end{rem}

\textbf{Proof.} a) We write $P(x, y) = Q(x, y) + L(x, y) + C =\ell$ where $Q(x, y) = a x^2 + b x y + cy^2$ with $a \not= 0$ is a totally positive quadratic form over $\mathcal{O}_F$ of discriminant $\Delta$,     $L(x, y) = \alpha_1 x + \alpha_2 y$ is a linear form over $\mathcal{O}_F$, and $C \in \mathcal{O}_F$.  
Let $\xi = (b \beta - 2 \alpha c)/\Delta$, $\eta = (b \alpha - 2 a \beta)/\Delta$. A little high school algebra shows
\begin{displaymath}
  P(x, y)   = \frac{(2a(x + \xi) + b(y+ \eta))^2 - \Delta (y+\eta)^2}{4a} + P(-\xi, -\eta). 
\end{displaymath}
Hence by a linear change of variables the equation $P(x, y) = \ell$ is equivalent to $$X^2 - \Delta Y^2 = 4a\Delta^2 (\ell - P(-\xi, - \eta)),$$ where $X$ and $Y$ satisfy certain congruence conditions modulo $4a|\Delta|^2$.  Clearly the norm of the left hand side is polynomial in $H(P )$ and $|\ell|_{\infty}$. We are now left with a norm form equation of the totally imaginary field $E = F(\sqrt{-|\Delta|})$ over the totally real field $F$, and the result follows. 

b) We use the representation \eqref{quasidiag} with $n=3$ together with the bounds \eqref{sizea}, \eqref{sizeh}, \eqref{allh}. By \eqref{lattice} we can choose $x_3$ in $\ll (\N(\ell)/\N (h_3)+1)^{1/2} \ll (\N(\ell)/\N(\Delta)^{1/3}+1)^{1/2}$ ways, and are left with an inhomogeneous  binary problem  for which part (a) applies.

c) This is proved in the same way. We choose $x_4$ and $x_3$ and are left with a binary problem. \qed \\

In the following lemma we denote by $\| . \|_2$ the usual Euclidean norm on $\Bbb{R}^n$, which is (in general) not induced by the inner product \eqref{inner}.  

\begin{lemma}\label{archlemma} Let $Q({\bf x}) = \frac{1}{2} {\bf x}^t A {\bf x}$ be a positive definite ternary quadratic form with real coefficients and eigenvalues $\gg 1$, let $ {\bf x} \in \Bbb{R}^3$ be such that $Q({\bf x}) = 1$ and let  $\ell > 0$,  $\eta > 0$.  

a) If $Q({\bf y}) = 1$ and $\langle {\bf y}, {\bf x}\rangle^2 \geq 1 - \eta$, then $\min_{\pm} \| {\bf y}\pm {\bf x}\|_2 \ll \eta^{1/2}$. 


b) Let $Q({\bf y}_i) = \ell$ for $i=1, 2, 3$ and assume  $|\langle{\bf y}_i, {\bf x}\rangle| \leq \ell^{1/2} \eta^{1/2}$.  Then $\det({\bf y}_1, {\bf y}_2, {\bf y}_3) \ll \ell^{3/2} \eta^{1/2}$. 

c) Let $Q({\bf y}) = \ell$, then $\| {\bf y} \|_2 \ll \ell^{1/2}$. \\
All implied constants are absolute (if the underlying number field is fixed). 
\end{lemma}   

\textbf{Proof.} The assertions are clear if $Q(\textbf{y}) = y_1^2 + y_2^3 + y_3^2$ and $\textbf{x} = (0, 0, 1)^t$ is the north pole. 
 In the general case, we write $\frac{1}{2}A = B^tB$ for some unique positive symmetric matrix $B \in \GL_3(\Bbb{R})$, so that $\|B\textbf{x}\|_2 = 1$. Let $S \in {\rm O}_3(\Bbb{R})$ be any orthogonal matrix with $SB\textbf{x} = (0, 0, 1)^t$. 
 Since $A$ has   eigenvalues $\gg 1$, the same holds for $B$ and hence for $SB$. For the proof of a) in the general case we conclude $\| \textbf{y} - \textbf{x} \|_2 \ll \| SB(\textbf{y} - \textbf{x})\|_2$, and for the two vectors $SB\textbf{y}$, $SB\textbf{x}$ the above special case applies. The other two parts are proved in the same way.\qed\\

The rather complicated proof of the next  lemma follows to some extent the argument in \cite[Lemma 2.1]{VdK}. It will be used for the bound in the eigenvalue aspect in Section \ref{eig-aspect}. Here and in the following we extend the norm $\N$ in an obvious way to a function $\Bbb{R}^d \rightarrow \Bbb{R}_{\geq 0}$. 

\begin{lemma}\label{VdKlemma} Let $Q( {\bf y}) = y_0^2 + \tilde{Q}(\tilde{ {\bf y}})$ with $\tilde{ {\bf y}} = (y_1, y_2, y_3)^t$ be an integral  positive definite quaternary quadratic form over $\mathcal{O}_F$ of discriminant $\Delta$, and let $\ell \in \mathcal{O}_F \setminus \{0\}$ be totally positive. Assume that $\tilde{Q}$ is reduced in the sense of Section \ref{reduction}. 
Let  $\eta = (\eta_1, \ldots, \eta_d) \in (0, 1]^d$. Let ${\bf x}_1, \ldots,   {\bf x}_d \in \Bbb{R}^3$ satisfy $\tilde{Q}^{\sigma_j}( {\bf x}_j) = 1$. Then the following two estimates hold:
\begin{equation}\label{bound1}
\begin{split}
 & \#\{  {\bf y} \in \mathcal{O}_F^4 \mid Q( {\bf y}) = \ell, \, (y^{\sigma_j}_0)^2 + \langle \tilde{ {\bf y}}^{\sigma_j},  {\bf x}_j\rangle^2 \leq \eta_j \ell^{\sigma_j} \} \\
 & \ll  \left( \N(\eta)^{1/2}\N(\ell)+1+ \min\left(\N(\ell)^{3/2}\N(\eta)^{1/2}, \N(\ell)^{1/2} \right)\right)(|\ell|_{\infty}\N \Delta)^{\varepsilon}. 
 \end{split} 
  \end{equation}
  and
  \begin{equation}\label{bound2}
  \begin{split}
 & \#\left\{ {\bf y} \in \mathcal{O}_F^4 \mid Q( {\bf y}) = \ell, \, \tilde{Q}^{\sigma_j}(\tilde{ {\bf y}}^{\sigma_j}) - \langle \tilde{ {\bf y}}^{\sigma_j},  {\bf x}_j \rangle^2  \leq \eta_j \ell^{\sigma_j} \right\} \\
 & \ll \left(1 + \min\left(\N(\eta)^{3/11} \N(\ell)^{12/11}  , \N(\ell) \right)\right)(|\ell|_{\infty}\N(\Delta))^{\varepsilon}. 
  \end{split}
   \end{equation}
   As before, all implied constant depend on $\varepsilon$ at most. 
\end{lemma}

\textbf{Proof.} We will frequently use Lemma \ref{archlemma} which is applicable because of \eqref{boundeig}. 

We start with the proof of \eqref{bound1}. By \eqref{lattice} there are  $\ll (\N(\eta)\, \N(\ell))^{1/2}$ choices for $y_0 \not=0$, and for each of them there are by Lemma \ref{lemma3}b at most $\ll \N(\ell)^{1/2}  (|\ell|_{\infty}\N \Delta)^{ \varepsilon}$ choices for $\tilde{\textbf{y}}$. This gives the first  term in \eqref{bound1}. 

We proceed to count the solutions with $y_0 = 0$. There are at most 2 linearly dependent solutions to $\tilde{Q}(\tilde{\textbf{y}}) = \ell$ (namely $\tilde{\textbf{y}}$ and $-\tilde{\textbf{y}}$), hence after adding 1 to the count of \eqref{bound1} we can assume   that there are at least two linearly independent solutions $\tilde{\textbf{y}}_1 = (y_{11}, y_{12}, y_{13})^t$, $\tilde{\textbf{y}}_2 = (y_{21}, y_{22}, y_{23})^t$, say, satisfying 
\begin{equation}\label{conditions}
  \tilde{Q}(\tilde{\textbf{y}}_\nu) = \ell, \quad \langle \tilde{\textbf{y}}^{\sigma_j}_\nu, \textbf{x}_j\rangle \leq (\eta_j \ell^{\sigma_j})^{1/2}, \quad \nu = 1, 2, \quad j = 1, \ldots, d.
\end{equation}  
Recall that by Lemma \ref{archlemma}c any solution $\textbf{y}$ to \eqref{conditions} satisfies $y_i^{\sigma} \ll (\ell^{\sigma})^{1/2}$ for $1 \leq i \leq 3$ and all $\sigma$.  Now any other  solution $\tilde{\textbf{y}}_3$ satisfies   $\det(\tilde{\textbf{y}}^{\sigma_j}_1, \tilde{\textbf{y}}^{\sigma_j}_2, \tilde{\textbf{y}}^{\sigma_j}_3) \ll (\ell^{\sigma_j})^{3/2} \eta_j^{1/2}$ by Lemma \ref{archlemma}b, as well as $\tilde{Q}(\tilde{\textbf{y}}_3)= \ell$. By \eqref{lattice} there are $1+\N(\ell)^{3/2} \N(\eta)^{1/2}$ choices for   the determinant (including 0). For a fixed value of the determinant and some $\tilde{\textbf{z}}_3 \in \mathcal{O}_F^3$ yielding this value, \emph{all} vectors yielding this determinant are of the form $\tilde{\textbf{y}}_3 = \tilde{\textbf{z}}_3 + a \tilde{\textbf{y}}_1 + b \tilde{\textbf{y}}_2 \in \mathcal{O}_F^3$ with $a, b \in K$.  Let $$\textbf{Y} := (Y_{23}, Y_{13}, Y_{12}) := \tilde{\textbf{y}}_1 \times \tilde{\textbf{y}}_2 \not= 0$$
be the cross product which is non-zero, since $\tilde{\textbf{y}}_1$  and $\tilde{\textbf{y}}_2$ are linearly independent. Assume without loss of generality that  $Y_{12} = y_{11}y_{22} - y_{12}y_{21} \not= 0$,   the other cases being analogous. Then $a\tilde{\textbf{y}}_1 + b\tilde{\textbf{y}}_2 \in \mathcal{O}_F^3$ implies that both $\alpha := Y_{12}a$ and $\beta := Y_{12} b$ are in $\mathcal{O}_F$, and $\tilde{Q}(Y_{12}\tilde{\textbf{y}}_3) = \tilde{Q}(Y_{12}\tilde{\textbf{z}}_3 + \alpha\tilde{\textbf{y}}_1 + \beta\tilde{\textbf{y}}_2 ) = Y_{12}^2\ell$ is a inhomogeneous binary problem in $\alpha, \beta$ with polynomial height in $|\ell|_{\infty}|$ and   the coefficients of $Q$ (recall that $Q$ is a reduced form), hence there are $\ll (|\ell|_{\infty}\N(\Delta))^{\varepsilon}$ solutions by Lemma \ref{lemma3}a. 
 
Alternatively, by Lemma \ref{lemma3}b we have the trivial bound $\ll \N(\ell)^{1/2} (|\ell|_{\infty}\N(\Delta))^{\varepsilon}$ for the number of solutions with $y_0 = 0$. Combining these two counts gives the last term in \eqref{bound1}. \\


We proceed to prove \eqref{bound2}.  Let $\delta_{1, j} \delta_{2, j} = \eta_j$.  We will fix $\delta_{1, j},\delta_{2, j}$ later and assume for the moment only $\N(\delta_1), \N(\delta_2) \ll \N(\ell)$.  It is enough to estimate the number of $\textbf{y}\in \mathcal{O}_F^4$ with $Q(\textbf{y}) = \ell$ satisfying $\tilde{Q}^{\sigma_j}(\tilde{\textbf{y}}^{\sigma_j}) \leq \delta_{1, j} \ell^{\sigma_j}$ and the number of $\textbf{y}\in \mathcal{O}_F^4$ with $Q(\textbf{y}) = \ell$ satisfying \begin{equation}\label{satisfying}(1 - \langle \tilde{\textbf{y}}^{\sigma_j}, \textbf{x}_j\rangle^2/\tilde{Q}^{\sigma_j}(\tilde{\textbf{y}}^{\sigma_j})) \leq \delta_{2, j}.\end{equation} 

We start with the latter. There are at most two solutions with $\tilde{\textbf{y}} = 0$. From now on we consider only solutions \textbf{y} with   $\tilde{\textbf{y}}\not= 0$. Among these we define an equivalence relation: we call $\textbf{y} = (y_0, \tilde{\textbf{y}}), \textbf{z} = (z_0, \tilde{\textbf{z}}) \in \mathcal{O}_F^4$ with $Q(\textbf{y}) = Q(\textbf{z}) = \ell$ equivalent if $\tilde{\textbf{z}} = c\tilde{\textbf{y}}$ for some $c \in K$. We claim that the cardinality of each equivalence class $[\textbf{y}]$ is small. This can be seen as follows:  Clearly  $(c) \subseteq (y_1, y_2, y_3)^{-1}$. Fix   a fractional ideal $\mathfrak{a} \supseteq \mathcal{O}_F$ with $\N \mathfrak{a} \ll 1$ in the  ideal class of $(y_1, y_2, y_3)$. Then $(c) \subseteq (y_1, y_2, y_3)^{-1} \mathfrak{a} = (\alpha)$, say, where $\N\alpha \asymp \N(y_1, y_2, y_3)^{-1} \gg \N(\ell)^{-1/2}$, and we choose   a generator such that  $|\alpha^{\sigma_j}| \asymp \N\alpha^{1/d}$, say.  Hence we can write $c = d \alpha$ with $d\in \mathcal{O}_F$. After multiplication with $1/\alpha$, the equation $Q(\textbf{z}) = \ell$ becomes an integral binary problem in $d$ and $z_0$ which by Lemma \ref{lemma3}a has $\ll (|\ell|_{\infty}\N(\Delta))^{\varepsilon}$ solutions. It is therefore enough to count the number of equivalence classes, and to this end we pick a set of representatives $\textbf{\textbf{y}}$; then by construction  the corresponding vectors $\tilde{\textbf{y}}$ are pairwise not collinear. 

From these representatives we select a vector $\textbf{y} \in \mathcal{O}_F^4$ such that 
\begin{displaymath}
  Y := \prod_{\sigma} \max(|(y_1)^{\sigma}|, |(y_2)^{\sigma}|)
\end{displaymath}
is maximal. We may assume $Y \not= 0$.  By Lemma \ref{archlemma}c we have $Y \ll \N(\ell)^{1/2}$. For any other (non-equivalent) solution $\textbf{z}  = (z_0, \ldots, z_3) = (z_0, \tilde{\textbf{z}}) \in \mathcal{O}_F^4$ we conclude from \eqref{satisfying} that
\begin{displaymath}
  \Bigl\langle \frac{\tilde{\textbf{y}}^{\sigma_j}}{\tilde{Q}^{\sigma_j}(\tilde{\textbf{y}}^{\sigma_j})^{1/2}}, \frac{\tilde{\textbf{z}}^{\sigma_j}}{\tilde{Q}^{\sigma_j}(\tilde{\textbf{z}}^{\sigma_j})^{1/2}}\Bigr\rangle \gg 1 - \delta_{2, j},
\end{displaymath}
and hence by Lemma 
\ref{archlemma}c  
\begin{equation}\label{det}
\begin{split}
  |(y_1z_2 - y_2z_1)^{\sigma_j}| &\leq \| \tilde{\textbf{y}}^{\sigma_j}  \|_2  \min_{\pm}  \| \tilde{\textbf{y}}^{\sigma_j} \pm \tilde{\textbf{z}}^{\sigma_j}  \|_2 \\
  &  =  \| \tilde{\textbf{y}}^{\sigma_j}  \|_2  (\ell^{\sigma_j})^{1/2}\min_{\pm}\Bigl\| \frac{\tilde{\textbf{y}}^{\sigma_j}}{\tilde{Q}^{\sigma_j}(\tilde{\textbf{y}}^{\sigma_j})^{1/2}} \pm \frac{\tilde{\textbf{z}}^{\sigma_j}}{\tilde{Q}^{\sigma_j}(\tilde{\textbf{z}}^{\sigma_j})^{1/2}} \Bigr \|_2\\
  &
      \ll \delta_{2, j}^{1/2} \| \tilde{\textbf{y}}^{\sigma_j} \|_2 (\ell^{\sigma_j})^{1/2}. 
   \end{split}   
\end{equation}
We first count solutions $\textbf{z}$ where $y_1z_2 - y_2z_1 \not= 0$. Then the non-zero principal ideal $(y_1z_2 - y_2 z_1)$ has norm 
\begin{displaymath}
  \ll D := \N(\delta_2)^{1/2} \prod_{\sigma} \| \tilde{\textbf{y}}^{\sigma} \|_2 \N(\ell)^{1/2} \ll \N(\delta_2)^{1/2} \N(\ell)
\end{displaymath}  
   and it is divisible by the ideal $\mathfrak{y} := (y_1, y_2)$. There are $\ll D/\N\mathfrak{y}$ such ideals. For each of these ideals the number of generators satisfying \eqref{det} is $\ll \N(\ell)^{\varepsilon}$ by \eqref{units}. Let us fix one such generator $g$ and count the number of $z_1, z_2 $ satisfying $y_1z_2 - y_2z_1 = g$ and $z^{\sigma}_1, z^{\sigma}_2 \ll (\ell^{\sigma})^{1/2}$ (by Lemma \ref{archlemma}c).  Assume there is a solution $(z_1, z_2)$. If $(\tilde{z}_1, \tilde{z}_2)$ is another solution, write $(z_1^{\ast}, z_2^{\ast}) := (\tilde{z}_1 - z_1, \tilde{z}_2 - z_2)$. Thus we need to count the number of $(z_1^{\ast}, z_2^{\ast})  \not= (0, 0)$ satisfying $(z_1^{\ast})^{\sigma}, (z_2^{\ast})^{\sigma} \ll (\ell^{\sigma})^{1/2}$ and $y_1z_2^{\ast} = y_2z_1^{\ast} $.  Fix an integral ideal $\mathfrak{a}$ in the ideal class of $\mathfrak{y}^{-1}$, and write $\mathfrak{y} \mathfrak{a} = (\alpha)$ with $\N\alpha \asymp \N\mathfrak{y}$. Then the fractional ideal $(y_2/\alpha)$ divides $(z_2^{\ast})$, and hence $z_2^{\ast} = y_2 w_2/\alpha$ for some $w_2 \in \mathcal{O}_F$. Similarly $z_1^{\ast} = y_1 w_1/\alpha$ for some $w_1 \in \mathcal{O}_F$. Inserting into the equation $y_1z_2^{\ast} = y_2z_1^{\ast}$ we see $w_1 = w_2 =: w\not=0$, say, and we have the bound $w^{\sigma} \ll (\ell^{\sigma})^{1/2} \alpha^{\sigma} \max(|y_1^{\sigma}|, |y_2^{\sigma}|)^{-1}$. By \eqref{lattice} the number of such $w$ is $\ll \N(\ell)^{1/2} \N\mathfrak{y}\, Y^{-1}$, and hence the number of $(z_1, z_2)$ for a fixed generator $g$  is $\ll 1 + N\ell^{1/2} \N\mathfrak{y} \, Y^{-1} \ll \N(\ell)^{1/2} \N\mathfrak{y} \, Y^{-1}$. Thus the total number of solutions $(z_1, z_2 )$ with $y_1z_2 - y_2z_1 \not= 0$ is
\begin{equation}\label{bounda}
   \ll \N(\ell)^{3/2+\varepsilon} \N(\delta_1)^{1/2} \, Y^{-1}.
\end{equation}     
Alternatively, by the definition of $\tilde{\textbf{y}}$ and $Y$, the principal ideals $(z_1)$ and $(z_2)$ have norm at most $Y$, and the number of generators whose conjugates are bounded by $O(\ell^{\sigma_j})$ is by \eqref{units} at most $\N(\ell)^{\varepsilon}$. Hence the number of $(z_1, z_2)$ is at most $\ll Y^2 \N(\ell)^{\varepsilon}$. Together with \eqref{bounda} we get the bound $\ll \N(\ell)^{1+\varepsilon} \N(\delta_2)^{1/3}$. Once $z_1, z_2$ are fixed, the equation $Q(\textbf{z}) = \ell$ is a binary problem in $z_0, z_3$, thus by Lemma \ref{lemma3}a our total count of (equivalence classes) $\textbf{z} \in \mathcal{O}_F^4$ with $y_1z_2 - y_2z_1 \not= 0$ is 
\begin{equation}\label{boundb}
  \ll \N(\delta_2)^{1/3} \N(\ell) \, (|\ell|_{\infty}\N \Delta)^{\varepsilon}.
\end{equation}
We now consider  the set  of all those (equivalence classes of) solutions $\textbf{z} \in \mathcal{O}_F^4$ with $y_1z_2 - y_2 z_1 = 0$. For any two such solutions $\textbf{z}, \textbf{w}$ we have by construction that $\tilde{\textbf{z}}, \tilde{\textbf{w}}$ are not collinear, hence $z_1z_3 - z_3z_1 \not = 0$. From the set of these vectors $\textbf{z}$ we select some $\textbf{z}^{\ast}$ such that
\begin{displaymath}
  Z := \prod_{\sigma} \max(|(z_1^{\ast})^{\sigma}|, |(z_3^{\ast})^{\sigma}|)
\end{displaymath}
is maximal. An argument identical to the above shows that the number of the (equivalence classes) $\textbf{z} \in \mathcal{O}_F^4$ such that $y_1z_3 - y_3 z_1 = 0$ satisfies the bound \eqref{boundb}.

It remains to count  the number of $\textbf{y}\in \mathcal{O}_F^4$ with $Q(\textbf{y}) = \ell$ satisfying $\tilde{Q}^{\sigma_j}(\tilde{\textbf{y}}^{\sigma_j}) \leq \delta_{1, j} \ell^{\sigma_j}$.   By Lemma \ref{lemma3}b the number of $\tilde{\textbf{y}}$ is 
\begin{equation}\label{boundc}
  \ll  1+ (\N(\delta_1) \N(\ell))^{3/2} (\max_{j}| \delta_{1, j}\ell^{\sigma_j}| \N(\Delta))^{\varepsilon}, 
\end{equation}  
and then $y_0$ is determined up to sign. Choosing 
  $\delta_{1, j}  = \eta_j^{2/11}(\ell^{\sigma_j})^{-3/11}$, the sum of \eqref{boundb} and \eqref{boundc} 
gives the first bound of \eqref{bound2}. The second bound follows trivially from Lemma \ref{lemma3}c.\qed\\

\section{Application of the pre-trace formula}\label{pretrace}
\subsection{The general set-up} As in \cite{BM1}, the proof of Theorem \ref{mainthm2} (which implies Theorem \ref{mainthm}) follows from an application of a pre-trace formula which we recall now. Our aim in this section is to state the pre-trace formula and to construct the amplifier. This will finally result in the important bound \eqref{term} to which we can apply the results from the preceding section. The pre-trace formula is   the spectral expansion of an automorphic kernel and works verbatim as in \cite[Section 5.1]{BM1} where the case $F = \Bbb{Q}$ is considered and to which we refer for more details. It features the matrix coefficients $p_{\bfm, \bfl}$ defined in Section \ref{archi_preliminaries}. Similar non-adelic treatments can be found in \cite{IS, VdK}.  The (only) new input compared to \cite{BM1}  is that we choose an amplifier that has support on ``balanced" algebraic integers  (i.e. algebraic integers whose various conjugates have roughly the same size) that generate principal prime ideals.   \\

Fix weights $\bfl \in \Bbb{Z}^d$ and a Hecke character $\chi$ as in Theorem \ref{mainthm2} (in particular $\chi$ is unramified at every place not dividing $\disc(\Oo)$). 
For  $\alpha=(\alpha_v)_v\in\BtAf$ such that $\alpha_v\in\Oo^\times_v\ \hbox{for $v\mid \disc(\Oo)$},$  let  $f_\alpha$ be the   function  supported  on
   $Z(\Af)\whOt \alpha\whOt$  that is   $1$ on $\whOt \alpha\whOt$ and  satisfies
  $$f_\alpha(\lambda h)=\ov\chi(\lambda)f_\alpha(h),\ \lambda\in Z(\Af).$$ 
 Let $\mcH(\Oo)$ be the convolution\footnote{for the convolution $f_1*f_2(g)=\int_{\PBAf}f_1(h)f_2(h^{-1}g)dh$,} (spherical) Hecke algebra generated by these bi-$\whOt$-invariant functions $f_{\alpha}$ on $\BtA$.  This algebra is commutative, and it follows from our assumptions (by $\whOt$-invariance) that $\vphi$ is an eigenfunction for the action of $\mcH(\Oo)$ by convolution: for $f\in\mcH(\Oo)$ one has
  $$f*\vphi(g)=\int_{\PB(\Af)}f(h)\vphi(gh)\ dh=\lambda_{\vphi}(f)\vphi(g).$$
  Given such an $f$, we consider more generally the convolution operator on $L^2(\Bt(F)\bash\BtA,\chi)$
  $$R(f):\psi(g)\mapsto R(f)(\psi)(g)=\int_\PBA f(h_f)p_{\bfm, \bfl}(h_\infty)\psi(gh)dh$$
  where we use as before the notation $h=h_fh_\infty\in\PBA$ and $p_{\bfm, \bfl}$ was defined in \eqref{matrix}.  
  This is an integral operator with kernel given by 
$$K_f(g,g')=\sum_{\gamma\in\PB(F)}f(g_f^{-1}\gamma g'_f)p_{\bfm, \bfl}(g_\infty^{-1}\gamma g'_\infty).$$
It decomposes into an orthonormal (finite) basis $\{\psi\}$ of $\whOt \rmK_{\infty}$-invariant $\mcH(\Oo)$-eigenfunctions of weight $\bfl$ with respect to $\rmK_\infty$ 
 containing $\vphi$, and  from  the normalization  of Haar measures  in Section \ref{measures} 
 one finds that
$$K_f(g,g')=\frac{1}{d_\bfm}\sum_{\psi}\lambda_{\psi}(f)\psi(g)\ov\psi(g').$$

Choosing $f$ appropriately, one can assume that $\lambda_{\psi}(f)\geq 0$ for any such $\psi$ and that  $\lambda_{\vphi}(f)$ is positive (and large): this is the principle of the amplification method. Taking $g=g'$, we obtain
$$|\vphi(g)|^2\leq \frac{|\bfm|}{\lambda_{\vphi}(f)}\sum_{\gamma\in\PB(F)}\bigl|f(g_f^{-1}\gamma g_f)p_{\bfm, \bfl}(g_\infty^{-1}\gamma g_\infty)\bigr|.$$
(Recall that $|\textbf{m}| \geq 1$ even for $\textbf{m} = 0$, which lightens our notation.) 

We now construct the amplifier $\lambda_{\vphi}(f)$ by a slight generalization of \cite[\S 5.2]{BM1} to the number field $F$ which takes into account  that the group of units $\mathcal{O}_F$ is in general infinite. Let $F_{\infty, +}^{\times}=\{(x_\sigma)_{\sigma\mid\infty}\in F^{\times}_{\infty},\ x_\sigma>0\}$ be the identity component of $F^\times_\infty$. We fix (once and for all) a fundamental domain $\mathcal{D}_0$ for the action of the totally positive units 
$U^{+}$ on the hyperboloid $\{y\in F_{\infty, +}^{\times} \mid\N y=1\}$. Let $F_{\infty, +}^{\text{diag}} = \{(y, \ldots, y) \in F_{\infty, +}^{\times}\}$. Then the cone $\mathcal{D} := F_{\infty, +}^{\text{diag}} \mathcal{D}_0$ is   a fundamental domain for the action of $U^{+}$ on $F_{\infty, +}^{\times}$. 

Given some parameter 
$L\geq 1$, consider the four sets 
\begin{displaymath}
\begin{split}
&  \mathcal{L}_1 := \{\ell \in \mathcal{O}_F\cap \mathcal{D} \mid \N(\ell) \in[L, 2L]\}, \\
& \mathcal{L}_2 :=  \{\ell \in \mathcal{O}_F \cap \mathcal{D} \mid \N(\ell) \in[L^2, (2L)^2]\},  \\
& \mathcal{L}_3 := \{\ell_1\ell_2^2 \mid L \leq \N(\ell_1), \N(\ell_2) \leq 2L, \, \ell_1, \ell_2 \in \mathcal{O}_F \cap \mathcal{D}\}, \\
&\mathcal{L}_4 := \{\ell^2 \mid L^2 \leq \N(\ell) \leq (2L)^2, \, \ell \in \mathcal{O}_F \cap \mathcal{D}\}. 
 \end{split} 
\end{displaymath}
where $\ell,\ell_1,\ell_2$ denote generators contained in $\mathcal{D}$ of integral \emph{principal   prime} ideals  $\mathfrak{p} \subseteq \mathcal{O}_F$ coprime with $\disc(\Oo)$. For any $\ell\in\mcL_{i}$ let 
$\alpha(\ell)=(\alpha(\ell)_v)_v\in\BtAf$ be  such that 
$\alpha(\ell)_v= \left(\begin{smallmatrix}\ell & 0 \\0 & 1\end{smallmatrix}\right)\in \GL_2(\mathcal{O}_{F,v})\simeq\Oo_v^\times$ if $v \nmid \disc(\Oo)$ and $\alpha(\ell)_{v}=1$ if $v \mid \disc(\Oo)$. Our selected function $f\in\mcH(\Oo)$ is then  the same linear combination of the $f_{\alpha(\ell)}$ as in \cite[p.\ 25]{BM1}. Precisely, for $r \in \mathcal{O}_F$ let
\begin{displaymath}
  c_r = \begin{cases} \text{sgn}( \lambda_{\vphi}(f_{\alpha( r)})), & r = \ell \in \mathcal{L}_1, \\
  \text{sgn}( \lambda_{\vphi}(f_{\alpha( r)})), & r = \ell^2, \ell \in \mathcal{L}_1, \\
  0, & \text{otherwise},
   \end{cases}
\end{displaymath}
let  $$\tilde{f} = \sum_{\ell}  c_{\ell} f_{\alpha(\ell )}$$
and define  $$f = \tilde{f} \ast \hat{\tilde{f}} = \sum_{d, \ell_1, \ell_2} c_{d\ell_1} c_{d\ell_2} f_{\alpha(\ell_1\ell_2)}.$$ 
For $g \in \PBA$ let  $\Oo_g$ be the order
$$\Oo_g:=\B(F)\cap g_f^{-1}\whO g_f$$
(where as before $g_f$ is the finite ad\`ele, the index does not refer to the function $f$). 
Observe that the  $\gamma\in\PB(F)$ satisfying $f_{\alpha(\ell)}(g_f^{-1}\gamma g_f)\not=0$  correspond to 
$\gamma\in\Oo_g$ with $\Nr(\tilde\gamma)=\ell$. This remark along with our choice of amplifier yields the bound
\begin{equation}\label{term}
  |\varphi(g)|^2 \ll |\bfm| L^{\varepsilon}\Bigl(\frac{1}{L} + \sum_{i=1}^4 \frac{1}{L^{2+i/2}} \sum_{ \ell \in \mathcal{L}_{i}}  \sum_{\substack{\gamma \in \Oo_g\\ \text{nr}(\gamma) = \ell}} | p_{\bfm, \bfl}(g_\infty^{-1}\gamma g_\infty) |\Bigr).
\end{equation}

Before we proceed, we make the important observation that  the quadratic form $Q_g$  associated to the order $\Oo_g$ is, in a suitable basis, of the form $y_0^2 + $ a ternary form in $y_1, y_2, y_3$, the latter corresponding to the restriction of the norm form to the traceless quaternions. 
  In particular, the assumptions of Lemma \ref{lemma2} and  Lemma \ref{VdKlemma} are satisfied. Moreover, by the discussion in Section \ref{quaternion}, we know that $\Oo_g$ is everywhere locally conjugate to the order $\Oo$, and hence $\Oo_g$ is an Eichler order and its associated quadratic form  $Q_g$ has the same discriminant $\Delta$ and the same level $\mfn$ as $Q=(\Oo,\nr)$. We set
\begin{displaymath}
  V:=\N(\disc\Oo_g)^{1/2}=\N(\Delta)^{1/2}=\N(\mfn), 
\end{displaymath}  
where the last equality follows from \eqref{level} and \eqref{disc}, and we note that $V = V_2^{1+o(1)}$ by \eqref{defV2}.

\subsection{Bound in the volume aspect}\label{volume} By the trivial bound $| p_{\bfm, \bfl}(g_\infty^{-1}\gamma g_\infty) |\leq 1$, we obtain
\begin{displaymath}
  |\varphi(g)|^2 \ll |\bfm| L^{\varepsilon}\Bigl(\frac{1}{L} + \sum_{i=1}^4  \frac{1}{L^{2+i/2}} \sum_{\ell \in \mathcal{L}_{i}} r_{Q_g}(\ell) \Bigr)
\end{displaymath}
where as before $Q_g$ is the   quaternary quadratic form associated with the order $\Oo_g$.  
We use Lemma \ref{lemma2} with $\N(\Delta)  = V^2$, $\N(\Delta)/\N(\mathfrak{n}) = V$. More precisely, we use \eqref{e3} for $i = 1$ resp.\ $i = 2$ with $y  \ll L$ resp.\ $y\ll  L^2$, we use \eqref{e2} for $i= 3$ with $y_1 = y_2 \ll L$ and we use \eqref{e1} for $i=4$ with $y \ll  L^2$.  In this way we obtain 
\begin{displaymath}
 |\varphi(g)|^2 \ll |\bfm|(VL)^{\varepsilon}\left(\frac{1}{L} + \frac{L^{1/2}}{V^{1/2}} + \frac{L^2}{V}\right).
 \end{displaymath}
 Choosing $L = V^{1/3}$ we find  
\begin{equation}\label{boundvol}
 \varphi(g) \ll_\eps |\bfm|^{1/2}V^{-1/6 + \varepsilon}.
\end{equation}

\subsection{Bound in the eigenvalue aspect}\label{eig-aspect} Let $x':=\rho_{g_\infty}(x^0) \in\B^0(F_\infty)\simeq (\Rr^3)^d$ be the $d$-tuple of unit vectors
$$x'_\sigma=\rho_{g_\sigma}(x^0_\sigma)=g_\sigma x^0_\sigma g_\sigma^{-1}$$
obtained by transforming the ``north pole"  $x_\infty^0$ by the rotation defined by $g_\infty$. The bi-$\rmK_\infty$-invariance of $ | p_{\bfm, \bfl}( . ) |$  implies that the function
$$\gamma\mapsto | p_{\bfm, \bfl}(g_\infty^{-1}\gamma g_\infty) |$$ depends on $\gamma$ (with $\nr(\gamma) = \ell \in \mathcal{O}_F$) only through the quantity
\begin{equation}\label{innerprod}
t = t(\gamma)=(t_\sigma)_{\sigma} \in [-1, 1]^d, \quad  t_\sigma=\peter{x'_\sigma,\gamma^\sigma x'_\sigma (\gamma^{\sigma})^{-1}}_{\B^0}= -1 + 2\frac{\langle \gamma^{\sigma}, 1\rangle_\B^2 + \langle \gamma^{\sigma}, x'_\sigma\rangle_\B^2}{\ell^{\sigma}} .
\end{equation}
The equality can be checked by brute force computation in the real quaternion algebra $\B(\Bbb{R})$, and by conjugation-invariance of the trace it is even enough to assume that $x_{\sigma}'$ is, say, the north pole in $\B^0(\Bbb{R})$. 

By  Lemma  \ref{decaylemma} we have
$$|p_{\bfm, \bfl}(g_\infty^{-1}\gamma g_\infty)|\ll_{\bfl} \prod_\sigma\min\Bigl(1,m_\sigma^{-1/2}(1-t_\sigma^2)^{-1/4}\Bigr).$$
We cut the $\gamma$-sum in \eqref{term} into dyadic (multi-dimensional) intervals according to 
\begin{displaymath}
\begin{split}
 & 1 \pm \langle x'_\sigma, \gamma^{\sigma} x_\sigma' (\gamma^{\sigma})^{-1}\rangle \asymp \eta_\sigma = 2^{-c_\sigma},  \quad 0 \leq c_\sigma <  C\\
 &  1 \pm \langle x'_\sigma, \gamma^{\sigma} x_\sigma' (\gamma^{\sigma})^{-1}\rangle \leq 2^{-C}
\end{split}  
\end{displaymath}
where $2^{-C} \asymp (LV|\bfm|)^{-100}$.  There are $O((LV|\bfm|)^{\varepsilon})$ such intervals.  
By \eqref{innerprod} the two conditions $1 \pm \langle x_\sigma', \gamma^{\sigma} x_{\sigma}' (\gamma^{\sigma})^{-1}\rangle \leq \eta_\sigma$ are equivalent to $\langle \gamma^{\sigma}, 1\rangle^2 + \langle \gamma^{\sigma}, x'_\sigma\rangle^2 \leq \eta_{\sigma} \ell^{\sigma}$ and (in the notation of Lemma \ref{VdKlemma}) $\tilde{Q}^{\sigma}(\tilde{ \gamma}^{\sigma}) - \langle \tilde{ \gamma}^{\sigma},  x_{\sigma}' \rangle^2  \leq \eta_{\sigma} \ell^{\sigma}$. Here we used again the special shape of our quadratic form. Hence by Lemma \ref{VdKlemma} we have 
\begin{displaymath}
\begin{split}
  \sum_{\substack{\gamma \in \Oo_g\\ \text{nr}(\gamma) = \ell \\ 1 \pm \langle x_\sigma', \gamma^{\sigma} x_{\sigma}' (\gamma^{\sigma})^{-1}\rangle \leq \eta_\sigma }}  1 \ll & \left(1+     \min\left( \N(\ell)^{1/2},  \N(\ell)^{3/2}\N(\eta)^{1/2}   \right) \right.\\
  &\left. + \min\left(\N(\eta)^{3/11}\N(\ell)^{12/11}, \N(\ell)   \right) \right) (V\N(\ell))^{\varepsilon}.  
 \end{split}
\end{displaymath}
Here we used  that $\N(\ell) \N \eta^{1/2} \leq \min\left(\N(\eta)^{3/11}\N(\ell)^{12/11}, \N(\ell)   \right)$.   If $ \eta_\sigma \leq 2^{-C}$ for at least one $\sigma$, then $\N(\eta) \ll (LV|\bfm|)^{-100}$, hence those $\gamma$ contribute  at most
\begin{equation}\label{eig1}
  \frac{|\bfm|}{L^{2+i/2}} L^{\min(i, 2)+\varepsilon} \ll  \frac{|\bfm|}{L^{1-\varepsilon}}. 
 \end{equation}
 to \eqref{term}. If $\eta_\sigma > 2^{-C}$ for all $\sigma$, then $\N(\eta) \gg (LV|\bfm|)^{-100d}$, and we are left with  
 \begin{equation}\label{opt}
\begin{split}
\ll & \sum_{i=1}^4 \frac{|\bfm| }{L^{2+i/2}} L^{\min(i, 2)} \max_{(LV|\bfm|)^{-100d} \ll \N(\eta) \leq 1}  \left(\frac{VL|\bfm|}{\N(\eta)}\right)^{\varepsilon}\min\left(1, \frac{1}{|\bfm|^{1/2}\N(\eta)^{1/4}}\right) \\
& \times  \left(1+     \min\left( L^{i/2},  L^{3i/2}\N(\eta)^{1/2}   \right)    + \min\left(\N(\eta)^{3/11}L^{12i/11},  L^{i} \right) \right).  
\end{split}
\end{equation}
We 
choose
\begin{displaymath}
  L := |\bfm|^{3/20}. 
\end{displaymath}
The discussion of \eqref{opt} as a function of $\eta$ with $0 \leq \N(\eta) \leq 1$ with this choice of $L$ is elementary, but a bit tedious.  The easiest way is to write $\N(\eta) = |\textbf{m}|^{\beta}$ with $\beta \leq 0$. Then \eqref{opt} becomes
\begin{displaymath}
\begin{split}
 (|\textbf{m}|V)^{\varepsilon}  &\sum_{i=1}^4 \max_{\beta \leq 0} \exp\Bigl\{(\log | \textbf{m}|) \Bigl(1 + \frac{3}{20}\Bigl(\min(i, 2) - 2 - \frac{i}{2}\Bigr)   + \min\Bigl(0, -\frac{1}{2}- \frac{1}{4}\beta\Bigr)\\
 & + \max\Bigl(0, \min\Bigl(\frac{i}{2} \cdot \frac{3}{20}, \frac{3i}{2} \cdot \frac{3}{20} + \frac{\beta}{2}\Bigr),  \min\Bigl(\frac{12i}{11}\cdot  \frac{3}{20} + \frac{3\beta}{11}, \frac{3i}{20}\Bigr)\Bigr)\Bigr)\Bigr\}.
\end{split}  
\end{displaymath} 
For each $i \in \{1, 2, 3, 4\}$, the exponent is a piecewise linear function in $\beta$ that is elementary to discuss. Figure 1 displays the four functions.  \\

\begin{center}\includegraphics[width=6cm]{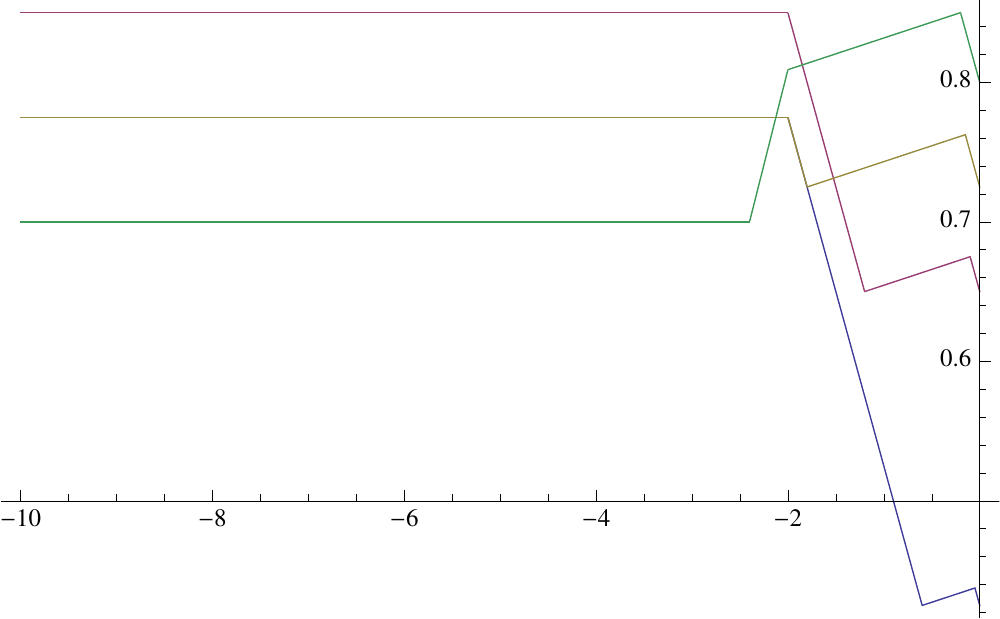}\end{center}

\noindent It is easy to check that the maximum is taken at $\eta \leq |\bfm|^{-2}$ and $i = 2$ and  at $\eta = |\bfm|^{-1/5}$ and $i = 4$ in which case \eqref{opt} is 
\begin{equation}\label{eig4}
  \ll |\bfm|^{17/20}(|\bfm|V)^{\varepsilon}. 
\end{equation}
 Combining \eqref{eig1} and \eqref{eig4} we obtain
\begin{equation}\label{finaleig}
  \varphi(g) \ll |\bfm|^{17/40}(|\bfm|V)^{\varepsilon}.   
\end{equation}

The bounds \eqref{boundvol} and \eqref{finaleig} prove \eqref{indiv} (observe that $3/80 < 1/27$). 

\subsection{Hybrid bound}

Combining \eqref{boundvol} and \eqref{finaleig} and writing $|\bfm| \asymp |{\lambda}|^{1/2}$, 
we obtain
\begin{displaymath}
\begin{split}
  \vphi(g) &\ll (|\lambda|V)^{\varepsilon} |\lambda|^{1/4} \times (V^{-1/6})^{9/29}( |\lambda|^{-3/80})^{20/29} \\
  &\ll |\lambda|^{1/4} (|\lambda|^{1/2}V)^{-3/58+\varepsilon } \ll  |\lambda|^{1/4} (|\lambda|^{1/2}V)^{-1/20 }. 
\end{split}
\end{displaymath}
This completes the proof of Theorem \ref{mainthm2}. 

\section{Bounds for automorphic forms on $3$-dimensional ellipsoids}\label{quaternary}

One of the main reasons for extending \cite{BM1} from $\Qq$ to general totally real number fields is that this will allow us to obtain bounds for automorphic forms associated to 
 {\em quaternary} quadratic spaces.  As was recalled in \S \ref{sec61} there is  a close relationship between automorphic functions associated to quaternary quadratic spaces and automorphic functions attached to ternary quadratic spaces but over a {\em quadratic} extension of the base field.
  
In the next section we define the arithmetic manifold $X^{(3)}(\Oo)$ as an adelic quotient of a special orthogonal group in four variables $G=\SO(\B')$ on which the automorphic forms considered in Theorem \ref{thmso4} live. In order to bound automorphic function on $G$, we pass to automorphic forms on the somewhat larger group $\widetilde{G}$, defined in \eqref{Gdef}.  Then the scene is prepared to copy the arguments from Section \ref{pretrace} and derive Theorem \ref{thmso4}.

\subsection{Automorphic forms associated to orthogonal groups in four variables}\label{autoSO4}

From now on we restrict to the following situation: $F$ is a fixed totally real number field, $E$ a fixed totally real quadratic extension of $F$ (possibly $F\times F$), $\B$
 is a totally definite quaternion algebra over $F$ and $\B'$ is the $4$-dimensional vector space of $\B\otimes_F E$ described in \S \ref{sec61}; as explained above, any totally definite quaternary quadratic space is similar to some $(\B',\Nr)$. In the following all implied constants may depend on $F$ and $E$. 
 
In the space  $(\B',\Nr)$ consider the following quadratic lattice: Let $\Oo\subset\B$ be an Eichler order, $\mathcal{O}_E$  a maximal order in $E$ (if $E=F\times F$ we take
$\mathcal{O}_E=\mathcal{O}_F\times\mathcal{O}_F$) and $\Oo_{\B_E}=\Oo\otimes_{\mathcal{O}_F}\mathcal{O}_E$; finally let $(\Lambda,\Nr)=(\Oo',\Nr)$ be the quadratic lattice with
$$\Oo'=\{z\in\Oo_{\B_E},\ {\sigma(z)^*}=z\}\subset\B'.$$
It is easy to see that $$\N(\disc\Oo')=\N(\disc\Oo)^{1+o_E(1)}.$$
We denote by $\rmK_{\Oo'}\subset \GAf$ the stabilizer of $\Oo'$ and by $\rmK_\infty=\prod_\sigma\rmK_\sigma\simeq\SO_3(\Rr)^d$ the stabilizer of the element $1\in\B'(F_\infty)$ in $\rmG(F_\infty)\simeq \SO_4(\Rr)^d$. We have the identification
\begin{equation}\label{defX3}
X^{(3)}(\Oo) := \rmG(F)\bash\rmG(\Aa)/\rmK_{\Oo'}\rmK_\infty \cong \bigsqcup_{i\in\mathrm{gen}(\Oo')} \Gamma_i\bash(S_3)^d,\ \Gamma_i\subset\SO_4(\Rr)^d
\end{equation}
where the disjoint union is indexed by the set of classes in the genus of $(\Oo',\Nr)$, and for a given representative of a genus class $(\Lambda_i,\Nr)$, $\Gamma_i$ is isomorphic to the subgroup of $\SO(\B')(F)$ which preserves the lattice $\Lambda_i$.  As above we are interested in a sup-norm bound of functions $\vphi$ on $X^{(3)}(\Oo)$ which are eigenfunctions of the Laplace operator 
$\Delta=(\Delta_\sigma)_\sigma$ and of a suitable algebra of Hecke operators $\mcH(\rmK_{\Oo'})$. The Laplace operator is normalized so that the eigenvalues have the shape 
\begin{equation}\label{eigenS3}
\lambda=(\lambda_\sigma)_\sigma,\ \hbox{where }\lambda_\sigma=-m_\sigma(m_\sigma+2), \hbox{ for $m_\sigma\in\Bbb{N}_{\geq 0}$.}
\end{equation}
These requirements imply that $\vphi$ is identified with a smooth, $\rmK_{\Oo'}\rmK_\infty$-invariant function, contained in the subspace $V_\pi\subset L^2(\rmG(F)\bash\GA)$ of an automorphic representation $\pi=\otimes_v\pi_v$ (as before we may assume that $\pi$ is infinite dimensional). 

By \cite[Thm.\ 4.13]{HS}, there exists an automorphic representation $\tpi\simeq \otimes_v\tpi_v$ of $\tGA$ (with unitary central character) such that $V_\pi\subset V^1_{\tpi|\GA}$ where $V^1_{\tpi|\GA}$ denotes the restriction to $\GA$ of a certain subspace of $V_\tpi$ (in most cases $V^1_{\tpi}=V_{\tpi}$). In other words, $\tpi$ is an automorphic representation on $\Bt(\Aa_E)$ whose central character on $Z_E(\Aa)\simeq\Aa_E^\times$ is trivial when restricted to $Z_F(\Aa)\simeq\At$. 

In particular, in order to bound $\vphi$ it is sufficient to bound $\tphi$, since obviously
$$\|\vphi\|_\infty\leq \|\tphi\|_\infty,$$
but we also need to compare the $L^2$-norms. We proceed to show that  $\|\vphi\|_2$ and $\|\tphi\|_2$ are of comparable size. First observe that since $\tphi$ is invariant by an open compact subgroup of $\tG(\Af)$ and by the subgroup $\{w\in\whO_{\B_E}^\times\mid \Nr(w)\in\At\}$, it is invariant under $\whO_{\B_E}^\times$. In particular $\tpi_v$ is unramified for every finite place not dividing $\N(\disc\Oo_{\B_E})$. It follows from this and \cite[Rem.\ 4.20]{HS} that for any choice for Haar measures
$$\frac{\peter{\vphi,\vphi}_\rmG}{\peter{\tphi,\tphi}_{\tG}}=C(\tpi)\frac{\vol(\rmG(F)\bash\GA)}{\vol(Z_E(\Aa)\tG(F)\bash\tGA)}$$
where $C(\tilde\pi)=\N(\disc\Oo_{\B_E})^{o(1)}$.
We fix the Haar measure compatibly with the fact that $\rmG$ is the kernel of the $F$-algebraic map 
\begin{eqnarray*}
\tG&\ra&  \mathbb{G}_{m,E}^1 = \{x \in \Bbb{G}_{m, E} \mid \N_{E/F}x = 1\}\\
w&\mapsto&\nr(w)/\sigma(\nr(w)),		
\end{eqnarray*}
by fixing a Haar measure on $\mathbb{G}_{m,E}^1(\Aa)$, and the Haar measure on $\PB(\Aa_E)$ so that $\PB(E_\infty)$ and the image of $\whO_{\B_E}^\times$ in $\PBAf$ have volume $1$. It follows from these choices that
\begin{displaymath}
\begin{split}
&\vol(Z_E(\Aa)\tG(F)\bash\tGA)=\NE(\disc\Oo_{\B_E})^{1/2+o_E(1)},\\
&\vol(\rmG(F)\bash\GA)=\vol(X^{(3)}(\Oo))=\N(\disc\Oo')^{1+o_E(1)}=\NE(\disc\Oo_{\B_E})^{1/2+o_E(1)},
\end{split}
\end{displaymath}
 and hence
\begin{equation}\label{comparenorms}
\frac{\peter{\vphi,\vphi}_\rmG}{\peter{\tphi,\tphi}_{\tG}}=\NE(\disc\Oo_{\B_E})^{o_E(1)}. 
\end{equation}

\subsection{The archimedean places} The algebra $E$ splits above the infinite places of $F$: $E_\sigma=\Rr\times\Rr$ and $\Bt_{E_\sigma}=\Bt_\sigma\times\Bt_\sigma$. By hypothesis $\tpi_\sigma$ admits for each $\sigma$ a non-zero
$\rmK_\sigma$-invariant vector. This implies that as a representation of $\Bu_\sigma\times\Bu_\sigma\simeq \SU\times\SU$ (which is a double cover of $\SO(\B'_\sigma)$)  we have 
 $$\tpi_\sigma\simeq\Pi_{m_{\sigma}}=\pi_{m_\sigma}\boxtimes \pi_{m_\sigma},$$ where $m_\sigma$ is defined by \refs{eigenS3}. For such a representation, the space of $\rmK_\sigma\simeq \SO_3(\Rr)$ invariant vectors is one-dimensional. Let $E_{m_\sigma}$ be such a non-zero vector, and for $g=(g_1,g_2)\in \Bu_\sigma\times\Bu_\sigma$ let
$$P_{m_\sigma}(g_\sigma)=\frac{\peter{g_\sigma.E_{m_\sigma},E_{m_\sigma}}}{\peter{E_{m_\sigma},E_{m_\sigma}}}$$
be the corresponding matrix coefficient. 
For $g=(g_\sigma)_\sigma\in \Bu(E_\infty)=\Bu(F_\infty)\times\Bu(F_\infty)$ let
$$P_{\bfm}(g)=\prod_\sigma P_{m_\sigma}(g_\sigma)$$ be the product of these matrix coefficients. 

\subsection{The pre-trace formula}
Let us (for notational simplicity) first assume that $E$ is a field. By the amplified pretrace formula for $\Bt(\Aa_E)$ we have for $g\in\Bt(\Aa_E)$ and a suitable parameter $L$ to be chosen in a moment 

\begin{equation}\label{term2}
\frac{|\tphi(g)|^2}{\peter{\tphi,\tphi}_{\tG}}\ll |\bfm|^2 L^{\varepsilon}\Bigl(\frac{1}{L} + \sum_{i=1}^4 \frac{1}{L^{2+i/2}} \sum_{ \ell \in \mathcal{L}_{i}}  \sum_{\substack{\gamma \in \Oo_{\B_E,g}\\ \text{nr}(\gamma) = \ell}} | P_{\bfm}(g_\infty^{-1}\gamma g_\infty) |\Bigr).
\end{equation}
This is the same expression as \eqref{term} except that the underlying field is now called $E$ rather than $F$ (this applies also to the definition of the sets $\mathcal{L}_i)$, and the matrix coefficient is different. In particular, $\Oo_{\B_E, g}$ is locally everywhere conjugated to $\Oo_{\B_E}$, and its level  
$\mfn_{\B_E}$ satisfies
$$\NE(\mfn_{\B_E}) \asymp_E \NE(\disc\Oo_{\B_E})^{1/2}$$
as in Section \ref{volume}. Let  $V  = \NE(\disc\Oo_{\B_E})^{1/2} =  \N(\disc\Oo')^{1+o(1)}$. 
Now the trivial bound $|P_{\textbf{m}}(g)| \leq 1$ together with the argument of Section \ref{volume} shows
$${|\tphi(g)|}\ll_{\eps} |\bfm| V^{-1/6+\eps}\|\tphi\|_2.$$
In view of \refs{eigenS3} and \refs{comparenorms} this may be rewritten
$${\|\vphi\|_\infty}\ll_{\eps}|\lambda|^{1/2} \vol(X^{(3)}(\Oo))^{-1/6+\eps}$$
if $\|\vphi\|_2 = 1$, which matches the generic bound   in the $\lambda$-aspect and improves it in the volume aspect. 

The split case $E= F \times F$ is very similar. In this case the sums over $\ell \in \mathcal{L}_i$ and over $\gamma$ in \eqref{term2} factor, and we can apply the argument of Section \ref{volume} for each factor. 
This completes the proof of Theorem \ref{thmso4}. 
 
\begin{rem} The function $P_m$ is bi-$\Delta\Bu_\sigma$ invariant (i.e.\  spherical) and may be expressed in terms of the character $\chi_m$ of the representation $\pi_m$, namely 
$$P_{m}(g)=P_{m}(g_1/g_2,1)=\chi_m(t)=\frac{1}{m+1}\frac{\sin((m+1)\theta)}{\sin(\theta)},$$
$$t=\cos(\theta)=\peter{g_1/g_2,1}_\B=\frac{1}2\tr(g_1/g_2).$$
In particular, one has
$$
P_{m}(g_1,g_2)\leq \min\Bigl(1,\frac1{(m+1)(1-t^2)^{1/2}}\Bigr). 
$$
Just as for Lemma \ref{decaylemma}, this decay property (as $t$ gets away from $\pm 1$) may be used together with an appropriate version of Lemma \ref{VdKlemma} to yield a hybrid bound for $\vphi$. We leave this to future work. 
\end{rem}

%


\begin{thebibliography}{EMN94}

\bib{BH}{article}{
   author={Blomer, Valentin},
   author={Holowinsky, Roman},
   title={Bounding sup-norms of cusp forms of large level},
   journal={Invent. Math.},
   volume={179},
   date={2010},
   number={3},
   pages={645--681},
   issn={0020-9910},
}
\bib{BM1}{article}{
 author={Blomer, V.},
  author={Michel, Ph.}
  title={Sup-norm bounds for eigenfunctions on ellipsoids},
 journal={IMRN},
 date={2011},
 volume={2011},
 number={21},
 pages={4934-4966},}

  

 \bib{Cas}{book}{
   author={Cassels, J. W. S.},
   title={Rational quadratic forms},
   series={London Mathematical Society Monographs},
   volume={13},
   publisher={Academic Press Inc. [Harcourt Brace Jovanovich Publishers]},
   place={London},
   date={1978},
   pages={xvi+413},
   isbn={0-12-163260-1},
}

\bib{DFIDuke}{article}{
   author={Duke, W.},
   author={Friedlander, J.},
   author={Iwaniec, H.},
   title={Class group $L$-functions},
   journal={Duke Math. J.},
   volume={79},
   date={1995},
   number={1},
   pages={1--56},
}

\bib{EMN}{article}{

   author={Erd{\'e}lyi, Tam{\'a}s},
   author={Magnus, Alphonse P.},
      author={Nevai, Paul},
   title={Generalized Jacobi weights, Christoffel functions, and Jacobi
   polynomials},
   journal={SIAM J. Math. Anal.},
   volume={25},
   date={1994},
   number={2},
   pages={602--614},
}

\bib{Fa}{book}{
   author={Faraut, Jacques},
   title={Analysis on Lie groups. An introduction},
   series={Cambridge Studies in Advanced Mathematics},
   volume={110},
   publisher={Cambridge University Press},
   place={Cambridge},
   date={2008},
   pages={x+302},
   isbn={978-0-521-71930-8},
}


\bib{GR}{book}{
   author={Gradshteyn, I. S.},
   author={Ryzhik, I. M.},
   title={Table of integrals, series, and products},
   edition={7},
   publisher={Elsevier/Academic Press, Amsterdam},
   date={2007},
   pages={xlviii+1171},
   isbn={978-0-12-373637-6},
   isbn={0-12-373637-4},
}


\bib{HT}{article}{
  author={Harcos, G.},
  author={Templier, N.},
  title={On the sup-norm of Maass cusp forms of large level. III},  
  journal={Math. Annalen, to appear},
 }




 \bib{HS}{article}{
   author={Hiraga, K.},
   author={Saito, H.}
   title={On $L$-packets for inner forms of $SL_{n}$},
  journal={Memoirs of the AMS},
   volume={215},
   date={2012},
}



\bib{Hum}{article}{
   author={Humbert, P.},
   title={Th\'eorie de la r\'eduction des formes quadratiques d\'efinies positives dans un corps alg\'ebrique $K$ fini},
   journal={ Comment. Math. Helv.},
   volume={12},
   date={1940},
   pages={263-306}
   }

\bib{IS}{article}{
   author={Iwaniec, H.},
   author={Sarnak, P.},
   title={$L^\infty$ norms of eigenfunctions of arithmetic surfaces},
   journal={Ann. of Math. (2)},
   volume={141},
   date={1995},
   number={2},
   pages={301--320},
   issn={0003-486X},
}

\bib{kitaoka}{book}{
    AUTHOR = {Kitaoka, Yoshiyuki},
     TITLE = {Arithmetic of quadratic forms},
    SERIES = {Cambridge Tracts in Mathematics},
    VOLUME = {106},
 PUBLISHER = {Cambridge University Press},
   ADDRESS = {Cambridge},
      YEAR = {1993},
     PAGES = {x+268},
     
}

 \bib{ponomarev}{article}{
   author={Ponomarev, Paul},
   title={Arithmetic of quaternary quadratic forms},
   journal={Acta Arith.},
   volume={29},
   date={1976},
   number={1},
   pages={1--48},
   issn={0065-1036},
}

\bib{Schur}{article}{
   author={Sarnak, Peter},
   title={Arithmetic quantum chaos},
   conference={
      title={The Schur lectures (1992) (Tel Aviv)},
   },
   book={
      series={Israel Math. Conf. Proc.},
      volume={8},
      publisher={Bar-Ilan Univ.},
      place={Ramat Gan},
   },
   date={1995},
   pages={183--236},
}

\bib{Sa}{article}{
   author={Sarnak, Peter},
   title={Letter to Morawetz},
   note={available at {\tt http://www.math.princeton.edu/sarnak}},
  }


 \bib{Te}{article}{
   author={Templier, N.},
   title={On the sup-norm of Maass cusp forms of 
large level},
   journal={Selecta Math.},
   volume={16},
   date={2010},
   pages={501--531},
}

 \bib{Tem1}{article}{
   author={Templier, N.},
   title={Large values of modular forms },
   journal={preprint {\tt arXiv:1207.6134}},
}

\bib{VdK}{article}{
  author={VanderKam, J. M. },
  title={$L^{\infty}$ norms and quantum unique ergodicity on the sphere},
  journal={IMRN},  
  volume={1997},
  date={1997},
  pages={329-347}
  }
  
  
  
  \bib{Ven}{article}{
   author={Venkatesh, Akshay},
   title={Sparse equidistribution problems, period bounds and subconvexity},
   journal={Ann. of Math. (2)},
   volume={172},
   date={2010},
   number={2},
   pages={989--1094},
}
  
  \bib{Vig}{book}{
   author={Vign{\'e}ras, Marie-France},
   title={Arithm\'etique des alg\`ebres de quaternions},
   language={French},
   series={Lecture Notes in Mathematics},
   volume={800},
   publisher={Springer},
   place={Berlin},
   date={1980},
   pages={vii+169},
   isbn={3-540-09983-2},
}

\bib{Vil}{book}{
   author={Vilenkin, N. Ja.},
   author={Klimyk, A. U.},
   title={Representation of Lie groups and special functions. Vol. 1},
   series={Mathematics and its Applications (Soviet Series)},
   volume={72},
   publisher={Kluwer Academic Publishers Group},
   place={Dordrecht},
   date={1991},
   pages={xxiv+608},
   isbn={0-7923-1466-2},
}



 \end{thebibliography}
\end{document}